\documentclass[11pt]{article}

\usepackage{a4wide}
\usepackage{enumerate}
\usepackage[greek, frenchb, english]{babel}
\usepackage[T1]{fontenc}        
\usepackage{mathrsfs} \usepackage{amssymb}
\usepackage{amsmath}
\usepackage{amsfonts}
\usepackage[dvips,pdftex
]{graphicx}
\usepackage{color} %\textcolor{orange{texte \`a colorier}
\usepackage[colorlinks=true, urlcolor=blue,linkcolor=blue, citecolor=blue]{hyperref}
\usepackage{aeguill}
\usepackage{rotating}
\usepackage{units}
\usepackage{graphicx}

\newtheorem{theo}{Theorem}[section]
\newtheorem{df}[theo]{Definition}%[section]
\newtheorem{lemme}[theo]{Lemma}%[section]
\newtheorem{Rem}[theo]{Remark}%[section]

\newtheorem{prop}[theo]{Proposition}
\newtheorem{hyp}{Hypothesis}

\newcommand{\CQFD}{\hfill $\square$}
\newcommand{\ind}{\mathbf{1}}

\def\real{\Bbb{R}}

\def\nit{\Bbb{N}}
\def\pp{\Bbb{P}}

\def\ee{\Bbb{E}}
\def\Bc{{\cal B}}

\def\Det{\mbox{Det}}
\def\Tr{\mbox{Tr}}

\title{Macroscopic analysis of determinantal random balls}

\author{
Jean-Christophe Breton\footnotemark[1], 
Adrien Clarenne\footnotemark[1] \
and
Renan Gobard\footnotemark[1]
}
\date{}
\begin{document}

\maketitle

\footnotetext[1]{IRMAR, Universit\'e de Rennes 1, 35042 Rennes cedex, France. Email: \{jean-christophe.breton, adrien.clarenne, renan.gobard\}@univ-rennes1.fr}

\begin{abstract}
We consider a collection of Euclidean random balls in $\real^d$ generated by a determinantal point process inducing interaction into the balls. 
We study this model at a macros\-copic level obtained by a zooming-out and three different regimes --Gaussian, Poissonian and stable-- are exhibited as in the Poissonian model without interaction. 
This shows that the macroscopic behaviour erases the interactions induced by the determinantal point process.
\end{abstract}

%%%%%%%%%%%%%%%%%%%%%%%%%%%%%%%%%%%%%%%%%%%%%%%%%%%%%%%%%%%%%%%%

\section*{Introduction}

A random balls model is a collection $\Bc$ of random Euclidean balls $B(x,r)=\{y\in\real^d : \|y-x\|\leq r\}$ whose centers $x\in\real^d$ and radii $r\in\real_+$ are generated by a stationary point process $N$ in $\real^d\times\real_+$. 
Such models are used to represent a variety of situation. 
Let mention a few of them.  
In dimension one, $\Bc$ can represent the traffic in a communication network. 
In this case, the (half-)balls are intervals $[x,x+r]$ and represent sessions of connection to the network, $x$ being the date of connection and $r$ the duration of connection.  
Such a model is investigated in \cite{MRRS2002} in a Poissonian setting, see also \cite{KT2008}. 
In dimension two, $\Bc$ can represent a wireless network with $x$ being the location of a base station emitting a signal with a range $r$ so that $B(x,r)$ represents the covering area of the station $x$ and the collection $\Bc$ gives the overall covering of the network, cf. \cite{YP2003}. 
The two-dimensional model is used also in imagery to represent Black and White pictures. 
In dimension three, such models are again used to represent porous media, for instance bones can be modeled in this way and an analysis of the model allows in this case to investigate anomalies such as osteoporosis, see \cite{BE2006}. 
Such random balls model is also known as grain-germ model with spherical grains in stochastic geometry, see the reference book \cite{CSKM2013}.

\medskip
In general in these models, one can think of at least two kinds of question, first, the geometrical --or morphological-- aspect of the collection $\Bc$ of balls and the corresponding continuum percolation problem; we refer to \cite{MR_perco1996} for this line of work. 
The second deals with the contribution of the model is some configuration and is the subject of this paper. 
For instance, the contribution in, say, a single site $y\in\real^d$, is given by the number of ball covering this site~$y$:
$$
\#\big\{B(x,r)\in \Bc :  y\in B(x,r)\big\}
=\sum_{B\in \Bc} \delta_y(B)
=\int_{\real^d\times\real_+} \delta_y\big(B(x,r)\big)\ N(dx,dr).
$$
Typically in the imagery setting ($d=2$), such a quantity gives the level of grey of pixel $y\in\real^2$, see \cite{BE2006}. 
More generally, one can consider the contribution of the model $\Bc$ into a configuration represented by a finite measure $\mu$, this contribution rewrites 
\begin{equation}
\label{eq:contribution}
M(\mu)=\int_{\real^d\times\real_+} \mu\big(B(x,r)\big)\ N(dx,dr).
\end{equation}
This is a shot-noise type functional and it will be the basic object of interest of this paper. 

\medskip
So far, these models have been investigated with a Poissonian generating mecanism, i.e. $N$ is a (homogeneous) Poisson point process (Ppp) with moreover center and radii behaviours being independent.
In addition to the above references, let mention \cite{KLNS2007} and \cite{BEK2010} where the $d$-dimensional model is investigated,  and \cite{BD2009} where weights are added to the balls.
A slight generalization is introduced in \cite{Gobard2014} where, still in a Poissonian paradigm, but non-homogeneous, the behaviours of the centers and of the radii are no more independent. 
Let also mention \cite{Heinrich_Schmidt1985}, \cite{Kluppelberg_Mikosch1995} and \cite{Lane1984} for asymptotics in related model for shot-noise processes 

\medskip
In the present paper, we go beyond the Poissonian setting and consider random balls generated by a stationary Determinantal process. 
As far as we know, except for the preliminary study \cite{Gobard_phd} where Ginibre point process (a special case of Determinantal process) is considered to generate the collection $\Bc$ and which is the very origin of this paper, this article presents the first study of a random balls model generated by a Determinantal point process, the so-called determinantal random balls. 
From a wireless network point of view, such a random mecanism is legitimate since it makes sense to install the stations not too close from one another. 
The repulsiveness of determinantal point processes justly realizes such a characteristic. 
From a modelling point of view, this choice has been recently explored in \cite{DZH}, \cite{MS2014} or \cite{LBDA2014}. 
In particular it is shown in \cite{DZH} that a thinned Ginibre point process is capable of modeling many of the actual cellular networks.

\medskip
Let us now be more specific about the macroscopic analysis provided in the sequel: we are interested in the behaviour of $M(\mu)$ in \eqref{eq:contribution} when a zoom-out in performed in the model. 
This zooming-out scheme offers at the limit a distant view of the model erasing the local specifities to make emerge only global characteristics. 
The scaling performed consists in $r\mapsto \rho r$ (with rate $\rho>0$) changing the ball $B(x,r)$ into $B(x,\rho r)$ and the zooming-out is performed with $\rho\to 0$. 
Denoting $M_\rho(\mu)$ for the contribution of the $\rho$-scaled model in $\mu$, a first-level description of $M_\rho(\mu)$ is given by its mean value 
$$
\ee\big[M_\rho(\mu)\big]
=\int_{\real^d\times\real_+} \mu\big(B(x,r)\big)\ n_\rho(dx,dr)
$$
where $n_\rho$ is the intensity measure of $N_\rho$, the image of $N$ by the scaling considered. 
A finer analysis is given by the fluctuations of $M_\rho(\mu)$ with respect to its mean value, i.e. the limit of 
\begin{equation}
\label{eq:Mrhomu1}
\frac{M_\rho(\mu)-\ee\big[M_\rho(\mu)\big]}{n(\rho)}
\end{equation}
for a proper normalization $n(\rho)$ when $\rho\to 0$. 
The limit above is investigated in distribution for each $\mu$, or, equivalently, because of the linear structure and thank to the Cram\'er-Wold device, in the finite-dimensional distributions (fdd) sense.  
Obviously, for the model not to vanish at the limit, the intensity, say $\lambda$, of the point process $N$ generating the balls has to be tuned accordingly. 
The relative behaviours of the scaling rate $\rho$ and of the balls intensity $\lambda$ will be responsible of the differents possible macroscopic regimes. 
A similar study has been done for the Poissonian random balls model, in which three different regimes -- Gaussian, Poissonian and stable-- appear at the limit, see \cite{KLNS2007, BEK2010}. 
Our study will justify that these regimes prevail for the determinantal random balls model, exhibiting thus a kind of robustness of theses regimes. 
Actually, since Poisson point processes are the universal limits of stationnary and ergodic point processes undergoing standard operations (independent thinning, dilatation), it is not surprising to recover similar asymptotics as the ones for the Poissonian model. 
We can even expect for these limits to be, in some way, universal. 

\medskip
The article is organized as follows. 
Section~\ref{sec:DetBall} gives a detailled presentation of the model investigated. 
The main results with the macroscopic behaviours (Theorems \ref{theo:largeBalls}, \ref{theo:intermediateBalls}, \ref{theo:smallBalls}) are stated and proved in Section~\ref{sec:asymptotics}. 
Examples are given in Section~\ref{sec:Exemples}.
Several final comments are gathered in Section~\ref{sec:comments} on zoom-in asymptotics, Ginibre point processes and $\alpha$-determinantal/permanental processes. 
Finally, Appendix~\ref{sec:mDpp} provides a very brief account on Determinantal point processes with the required results for our analysis.

%%%%%%%%%%%%%%%%%%%%%%%%%%%%%%%%%%%%%%%%%%%%%%%%%%%%%%%%%%%%%%%%%%%%

\section{Determinantal random balls model}
\label{sec:DetBall}

The model considered is a collection ${\Bc}$ of random (Euclidean) balls $B(x,r)=\big\{y\in\real^d : \|y-x\|\leq r\big\}$ whose centers $x\in\real^d$ and radii $r\in\real_+$ are generated by a marked stationary determinantal point process (Dpp) $\Phi$ on $\real^d\times \real_+$. 
In this section, we describe thoroughly the model and we refer to the Appendix~\ref{sec:mDpp} for more details on Dpp, in particular see the definition in Def.~\ref{def:Dpp}. 
First, consider a stationary Dpp $\phi$ with a kernel $K$ with respect to the Lebesgue measure $Leb$ satisfying $K(x,y)=K(x-y)$ (for simplicity, we use the same letter $K$ for two different functions), moreover we assume that the map ${\bf K}$ given for all $x\in \real^d$ and any $f\in L^2(\real^d, dx)$
\begin{equation}
\label{eq:KK}
{\bf K}f(x)=\int_{\real^d}K(x,y)f(y)\ dy
\end{equation}
satisfies the following hypothesis
\begin{hyp} 
\label{hyp:det1}
The map ${\bf K}$ in \eqref{eq:KK} is a bounded symmetric integral operator ${\bf K}$ from $L^2(\real^d,dx)$ into $L^2(\real^d,dx)$ with spectrum included in $[0,1[$. 
Moreover, ${\bf K}$ is locally trace-class, i.e. for all compact $\Lambda \subset E$, the restriction ${\bf K}_\Lambda$ of ${\bf K}$ on $L^2(\Lambda,\lambda)$ is of trace-class. 
\end{hyp}
This point process $\phi$ generates the centers of the balls and as a Dpp exhibits repulsiveness between its particles. 
To obtain balls, attach to each center $x$ a (positive) mark interpreted as a radius $r$ and independently and identically distributed according to $F$, assumed to admit a probability density $f$.  
The collection of these marks and of the Dpp $\phi$ forms a marked Dpp~$\Phi$. 
According to Proposition~\ref{prop:MDPP}, $\Phi$ is a Dpp on $\real^d\times \real_+$ with kernel
$$
\widehat{K}\big((x,r),(y,s)\big)= \sqrt{f(r)}K(x,y)\sqrt{f(s)},
$$
with respect to the Lebesgue measure.
Denoting $\Phi$ as well for the marked Dpp $s\{(X_i, R_i)\}$ as for the associated random measure $\sum_{(X,R)\in \Phi}\delta_{(X,R)}$, for any point $y\in\real^d$, the number of balls containing $y$ is given by:
$$
M(y)=\sum_{(X_i,R_i)\in \Phi} \ind_{B(X_i,R_i)}(y)=\int_{\real^d\times\real_+}\delta_y\big(B(x,r)\big)\ \Phi(dx,dr),
$$
where, for any set $A$, $\delta_y(A)=\ind_A(y)$. 
Identifying $y$ with the Dirac measure $\delta_y$, the previous definition actually extends from any Dirac measure $\delta_y$ to any suitable (signed) measure $\mu$ on  $\real^d$, defining the contribution of the model in such a configuration $\mu$ as the following field 
\begin{equation}
\label{eq:procM}
M(\mu)=\int_{\real^d\times\real_+}\mu\big(B(x,r)\big)\ \Phi(dx,dr).
\end{equation}
Note that from a mathematical point of view, it is not required for the measure $\mu$ to be positive and signed measures can be considered. 
However, in order to ensure that $M(\mu)$ in \eqref{eq:procM} is well defined, we restrain to measures $\mu$ with finite total variation  (see below Proposition~\ref{prop:Mesp}). 
In the sequel, we note ${\cal Z}(\real^d)$ the set of signed measure $\mu$ on $\real^d$ with finite total variation $\|\mu\|_{var}(\real^d)<+\infty$.
Moreover as in \cite{KLNS2007}, assume the following assumption on the radius behaviour, for $d<\beta<2d$, 
\begin{equation}
\label{eq:tails}
f(r)\underset{r\to +\infty}{\sim} \frac{C_\beta}{r^{\beta+1}}, \qquad
r^{\beta+1} f(r)\leq C_0.
\end{equation}
Since $\beta>d$, condition \eqref{eq:tails} implies that the mean volume of the random ball is finite:
\begin{equation}
\label{eq:cfvol}
v_d\int_0^{+\infty}r^df(r)\ dr<+\infty.
\end{equation}
where $v_d=Leb\big(B(0,1)\big)=\pi^{d/2}/\Gamma(d/2+1)$ is the Lebesgue measure of the unit ball of $\real^d$.
On the contrary, $\beta<2d$ implies that $F$ does not admit a moment of order $2d$ and the volume of the balls has an infinite variance. 
This is responsible of some kind of long-range dependence in the model, see \cite[p.~530]{KLNS2007} and is in line with communication network models which exhibit interference.
The asymptotics condition in \eqref{eq:tails} is of constant use in the following.  

\begin{prop}
\label{prop:Mesp}
Assume \eqref{eq:tails} is in force. 
For all $\mu \in {\cal Z}(\real^d)$, $\ee\left[\vert M(\mu) \vert \right]<+ \infty$. 
As a consequence, $M(\mu)$ in \eqref{eq:procM} is almost surely well defined for all $\mu\in{\cal Z}(\real^d)$.
\end{prop}
\begin{Proof}
Using properties of functionals of random measures (see Section~9.5 in \cite{DVJ1}), we have:
$$
\ee\big[ |M(\mu)| \big]
=\int_{\real^d\times\real_+} \big| \mu\big(B(x,r)\big)\big| \widehat{K}\big((x,r),(x,r)\big)\ dxdr.
$$
Since $\widehat{K}\big((x,r),(x,r)\big)=K(0)f(r)$, writing 
$\mu\big(B(x,r)\big)=\int_{\real^d} \ind_{B(y,r)}(x)\ \mu(dy)$,
we have
\begin{eqnarray*}
\ee\big[ |M(\mu)| \big]
&\leq& \int_{\real^d\times\real_+} \int_{\real^d} \ind_{B(y,r)}(x)\ |\mu|(dy)\ K(0)f(r)\ dxdr \\
&\leq& K(0)
\int_{\real^d}\int_0^{+\infty} \Big(\int_{\real^d} \ind_{B(y,r)}(x)\ dx\Big) f(r)\ dr\ |\mu|(dy) \\
&\leq& K(0)
 Leb\big(B(0,1)\big) \int_0^{+\infty} r^df(r)\ dr\ \int_{\real^d}|\mu|(dy) \\
&\leq& v_d\|\mu\|_{var}K(0) \Big(\int_0^{+\infty} r^df(r)\ dr\Big). 
\end{eqnarray*}
This concludes the proof thanks to conditions \eqref{eq:cfvol}, due to \eqref{eq:tails}.
\CQFD
\end{Proof}

%%%%%%%%%%%%%%%%%%%%%%%%%%%%%%%%%%%%%%%%%%%%%%%%%%%%

\section{Asymptotics}
\label{sec:asymptotics}

We now detail our zooming-out procedure. 
This procedure acts accordingly both on the centers and on the radii (equivalently on the volume of the balls). 
First, a scaling $S_\rho: r\mapsto \rho r$ of rate $\rho<1$ changes balls $B(x,r)$ into $B(x,\rho r)$; this scaling changes the distribution $F$ of the radius into $F_\rho=F\circ S_\rho^{-1}$. 
Second, the intensity of the center is simultaneously adapted; to do this, introduce a scaled version $\phi_\rho$ of the Dpp $\phi$, with kernel $K_\rho$ with respect to the Lebesgue measure. 
In order to be in line with previous model ball, we introduce $\lambda(\rho)$ given by 
\begin{equation} 
\label{eq:lambdarho}
K_\rho(0)=\lambda (\rho)K(0)
\end{equation}
with $\lim_{\rho\to 0}\lambda(\rho)=+\infty$.
For technical purpose, suppose also that, for any $\rho>0$,
\begin{equation}
\label{eq:controlintKrho} 
\int_{\real^d}|K_\rho(x)|^2\ dx \underset{\rho \rightarrow 0}{=} \mathcal{O}\big(\lambda(\rho)\big).
\end{equation}
\begin{Rem}
{\rm
The quantity $\lambda (\rho)$ introduced in \eqref{eq:lambdarho} can be interpreted as the intensity of the balls. 
Then $\lambda(\rho)\to+\infty$ indicates there is more and more balls when zooming-out ($\rho\to 0$). 
The condition \eqref{eq:controlintKrho} gives a control of $K_\rho(x)$ for $x\not=0$ and, roughly speaking, means that the correlation of the balls is controled by the intensity of the balls.  
}
\end{Rem}

\medskip\noindent
In summary, the zoom-out procedure consists in considering a new marked Dpp $\Phi_\rho$ on $\real^d\times \real_+$ with kernel:
$$
\widehat{K}_\rho\big((x,r),(y,s)\big)
=\sqrt{\frac{f(r/\rho)}{\rho}}K_\rho(x,y)\sqrt{\frac{f(s/\rho)}{\rho}},
$$
with respect to the Lebesgue measure. 
The scaled version of $M(\mu)$ is then the field
$$
M_{\rho}(\mu)=\int_{\real^d\times\real_+} \mu\big(B(x,r)\big)\ \Phi_\rho(dx,dr).
$$
In the sequel, we are interested in the fluctuations of $M_{\rho}(\mu)$ with respect to its expectation 
$$
\ee\big[M_{\rho}(\mu)\big]=\int_{\real^d\times\real_+}\mu\big(B(x,r)\big)\ K_\rho(0)\frac{f(r/\rho)}{\rho}\ dxdr
$$
and we introduce
\begin{equation}
\label{eq:procMM}
\widetilde M_\rho(\mu)
=M_\rho(\mu)-\ee\big[M_\rho(\mu)\big]
=\int_{\real^d\times\real_+} \mu\big(B(x,r)\big)\ \widetilde\Phi_\rho(dx,dr),
\end{equation}
 where $\widetilde\Phi_\rho$ stands for the compensated random measure associated to $\Phi_\rho$. 
 
%%%%%%%%%%%%%%%%%%%%%%%%%%%%%%%%%%%%%%%%%%%%

\subsection*{Heuristics}

The asymptotic behavior of $\widetilde M_\rho(\mu)$ when $\rho \to 0$ depends on how the scaling rate $\rho $ and the intensity $\lambda(\rho)$ are tuned. 
Roughly speaking, three regimes appear according to 
$\rho\to 0$ faster, slower or well-balanced with respect to $\lambda(\rho)\to+\infty$. 
Heuristically, the key quantity ruling these regimes is the mean number of large balls, say balls of radii larger than $1$ and, say, containing $0$:
 \begin{eqnarray*}
\lefteqn{\ee\Big[\#\big\{(x,r) \in \Phi_\rho\ : \ 0\in B(x,r), r>1 \big\} \Big]
}\\
&=&\int_{\{(x,r)\ :\ 0\in B(x,r), r>1\}} \widehat{K}_\rho\big((x,r),(x,r)\big)\ dxdr
=\int_1^{+\infty} \int_{B(0,r)} K_\rho(x,x)\frac{f(r/\rho)}{\rho}\ dxdr 
\\
&=&\int_{1/\rho}^{+\infty} \int_{B(0,\rho u)} \lambda(\rho)K(0)\ dx \ f(u)\ du 
\sim v_d K(0)\lambda(\rho)\rho^d \int_{1/\rho}^{+\infty} u^{-1-\beta+d}\ du 
\\
&\sim&\frac{v_d K(0)}{\beta-d}\lambda(\rho)\rho^\beta 
 \end{eqnarray*}
using \eqref{eq:lambdarho}, \eqref{eq:tails}. 
Thus the balance between $\rho\to 0$ and $\lambda(\rho)\to+\infty$ is ruled by $\lambda(\rho)\rho^\beta$ and the three scaling regimes are the following when $\rho \to 0$:
\begin{itemize}

\item Large-balls scaling: $\lambda(\rho)\rho^{\beta} \to +\infty$. 
Roughly speaking, large balls prevail at the limit and they shape the limit according to some kind of CLT. 
Moreover, since the large balls overlap, this regime yields dependence at the limit. 
In other words, the limit $\lambda(\rho)\rho^\beta\to+\infty$ acts as if $\lambda\to+\infty$ first and $\rho\to 0$ next ; the first limit ($\lambda\to+\infty$) corresponds to the superposition of a large number of (overlapping) balls, which in line with a CLT argument, produces a Gaussian limit (with dependence), the second limit ($\rho\to 0$) only shapes the covariance of the Gaussian field.
In this context, the proper normalization will be $n(\rho)=\sqrt{\lambda(\rho)\rho^{\beta}}$. 
See Section~\ref{sec:large_scaling}.

\item Intermediate scaling: $\lambda(\rho)\rho^{\beta} \to a \in ]0,+\infty[$. 
Roughly speaking, there is a proper balance between large and small balls and somehow the limit is incompletely taken and it only consists in an alteration of the generating point process with a dissolving of the interaction resulting in a Poisson point process. 
In this context, the proper normalization will just be a constant. 
See Section~\ref{sec:intermediate_scaling}.

\item Small-balls scaling: $\lambda(\rho)\rho^{\beta} \to 0$. 
Roughly speaking, small balls prevail. 
In other words the limit $\lambda(\rho)\rho^\beta\to 0$ acts as if $\rho\to 0$ first and $\lambda\to+\infty$ next. 
The first limit $\rho\to 0$ is a scaling killing the overlapping and thus producing independence at the limit. 
Next, with the second limit ($\lambda\to+\infty$) the heavy-tails of $F$ enter the picture: the contribution of the non-overlapping balls are in the domain of attraction of a stable distribution producing a stable regime. 
Moreover, the index of stability $\gamma$ can be heuristically derived as follows: for a smooth measure $\mu$, we have $\mu\big(B(x,r)\big)\asymp cr^d$ with $(\beta/d)$-regular tails under \eqref{eq:tails} responsible of the stability of index $\gamma=\beta/d$.
See Section~\ref{sec:small_scaling}.

\end{itemize}
 
%%%%%%%%%%%%%%%%%%%%%%%%%%%%%%%%%%%%

\subsection*{General strategy}
\label{sec:strategy}

For the three regimes, the proofs will follow the same idea in Sections~\ref{sec:large_scaling}, \ref{sec:intermediate_scaling}, and  \ref{sec:small_scaling} below, and the general strategy is presented. 
 The main tool to study the so-called determinantal integral  \eqref{eq:procM} or \eqref{eq:procMM} (integral with respect to a determinantal random measure) is the Laplace transform given in  Theorem~\ref{theo:Laplace_Dpp}. 
However, this result applies for compactly supported integrands which is not the case in our case with $(x,r)\mapsto \mu\big(B(x,r)\big)$ (when $r\to+\infty$, $\mu\big(B(x,r)\big)\to\mu(\real^d)$).
As a consequence, we consider the following auxiliary truncated process: 
\begin{equation}
\label{eq:procMR}
M_\rho^R(\mu)=\int_{\real^d\times\real_+} \mu\big(B(x,r)\big) \ind_{\{ r\leq R\}}\ \Phi_\rho(dx,dr),
\end{equation}
and the associated compensated determinantal integral $\widetilde M_\rho^R(\mu)$. 
Then for a positive compactly supported measure $\mu$, $(x,r)\mapsto \mu\big(B(x,r)\big )\ind_{\{r\leq R\}}$ is indeed a compactly supported function. 
In the following, we thus restrain ${\cal Z}(\real^d)$ to ${\cal Z}_c^+(\real^d)$ the set of positive compactly supported measures on $\real^d$ with finite total variation. 
The relevance in introducing this auxiliary process appears in the following result:
\begin{prop}
\label{prop:cvmR}
Assume \eqref{eq:tails} and \eqref{eq:lambdarho}. 
For all $ \mu \in \mathcal{Z}_c^+(\real^d)$ and for all $\rho>0$, $M_\rho^R(\mu)$ converges in $L^1$ when $R\to +\infty$ to $M_\rho(\mu)$. 
Moreover, in the intermediate and the small-balls scalings, there exists a constant $\rho_1>0$, independent of $R$, such that this convergence is uniform in $\rho$ for $\rho \in (0, \rho_1)$.
\end{prop}
\begin{Proof}
Let $\mu \in {\cal Z}_c^+(\real^d)$. 
By the monotone convergence theorem $M_\rho^R(\mu)\nearrow M_\rho(\mu)$ when $R\to+\infty$ and by the dominated convergence theorem $M_\rho^R(\mu)\to M_\rho(\mu)$ in $L^1$. 
Next, we have 
$$
M_{\rho}(\mu)-M_{\rho}^R(\mu)
=\int_{\real^d\times\real_+}\mu\big(B(x,r)\big) \ind_{\{r> R\}}\ \Phi_\rho(dx,dr),
$$
and
\begin{eqnarray*}
\ee\Big[\Big|M_{\rho}(\mu)-M_{\rho}^R(\mu)\Big|\Big]
&=&
\ee\Big[\int_{\real^d\times\real_+}\mu\big(B(x,r)\big)\ind_{\{r> R\}}\ \Phi_\rho(dx,dr)\Big]
\\
&=&\int_{\real^d\times\real_+}\mu\big(B(x,r)\big)\ind_{\{r> R\}} \widehat K_\rho\big((x,r),(x,r)\big)\ dxdr\\
&=&\int_{\real^d}\int_R^{+\infty} \mu\big(B(x,r)\big)  K_\rho(x,x)\frac{f(r/\rho)}{\rho}\ dxdr\\
&=&  \lambda(\rho)K(0)\int_{\real^d}\int_R^{+\infty}
\mu\big(B(x,r)\big)\frac{f(r/\rho)}{\rho}\ dxdr
\end{eqnarray*}
But with Fubini theorem and a change of variables
\begin{eqnarray*}  
\int_{\real^d}\int_R^{+\infty}\int_{\real^d} \ind_{B(x,r)}(y) \frac{f(r/\rho)}{\rho}\ \mu(dy)dxdr
&=&  \int_{\real^d}\int_R^{+\infty}v_dr^d\frac{f(r/\rho)}{\rho}\ dr\mu(dy)\\
&=& v_d\mu(\real^d)\rho^d\int_{R/\rho}^{+\infty}u^d f(u)\ du.
\end{eqnarray*}
From \eqref{eq:tails}, 
we have $f(u)\leq C_0/u^{\beta+1}$ and when $\rho<1$, 
$$
\rho^d\int_{R/\rho}^{+\infty}u^d f(u)\ du\leq \rho^d \int_{R/\rho}^{+\infty}u^d \frac{C_0}{u^{1+\beta}}\ du 
=\frac{C_0}{\beta-d}R^{d-\beta} \rho^\beta
$$
so that
$$
\ee\Big[\Big|M_{\rho}(\mu)-M_{\rho}^R(\mu)\Big|\Big]
\leq \frac{C_0}{\beta-d}R^{d-\beta}\lambda(\rho)\rho^\beta K(0)v_d\mu(\real^d)
$$
which goes to $0$ uniformly in $\rho\in]0,1[$ under the intermediate and small scaling since $\lambda(\rho)\rho^\beta$ is bounded in these case. 
\CQFD
\end{Proof}

\medskip\noindent
This uniform convergence is crucial in order to interchange the limit in $\rho$ and the limit in $R$ whenever $\lim_{\rho\to 0} \widetilde{M}_\rho^R(\mu)$ exists:
\begin{equation}
\label{eq:intervertion}
\lim_{\rho \to 0} {\cal L}\big(\widetilde{M}_\rho(\mu)\big)
=\lim_{\rho \to 0}\lim_{R\to +\infty}  {\cal L}\big(\widetilde{M}_\rho^R(\mu)\big)
=\lim_{R\to +\infty} \lim_{\rho \to 0} {\cal L}\big(\widetilde{M}_\rho^R(\mu)\big).
\end{equation}

\medskip\noindent
The strategy is now clear to obtain $\lim_{\rho \to 0} \widetilde{M}_\rho(\mu)$: (i) first, take $\lim_{\rho \to 0} \widetilde{M}_\rho^R(\mu)$ and (ii) next take the limit in $R\to+\infty$. 
In order to realize (i), use the Laplace transform of Dpp \eqref{eq:Laplace_Dpp} and the expansion \eqref{eq:Fredholm_det} of the corresponding Fredholm determinant. 
In this expansion, the first term (for $n=1$) is identified as a Poissonian term for which the asymptotics of the Poissonian model applies and the remaining terms ($n\geq 2$) are shown to be asymptotically negligible. 
Next, (ii) properly shapes the limit with $R\to+\infty$.

\medskip\noindent
However in order to realize (i), it is required to investigate the convergence of $\widetilde{M}_\rho^R(\mu)$ when $\rho \to 0$ on a restricted class of measures $\mu$ that we introduce now.
\begin{df}
\label{def:Mbeta}
The set ${\cal M}_\beta^+$ consists of positive measures $\mu \in {\cal Z}_c^+(\real^d)$ such that there exist two real numbers $p$ and $q$ with $0<p<\beta< q\leq 2d$ and a positive constant $C_\mu$ such that 
\begin{equation}
\label{eq:Mbeta}
\int_{\real^d}\mu\big(B(x,r)\big)^2\ dx\leq C_\mu\big( r^{p}\wedge r^{q}\big),
\end{equation}
where $a\wedge b=\min(a,b)$.
\end{df}
The control in \eqref{eq:Mbeta} by both $r^p$ and $r^q$ is required to ensures our quantity are well defined (see Proposition~\ref{prop:Mcal}-(i)); however in the sequel, only the control by $r^q$ will be used.
This definition is reminiscent of ${\cal M}_{2,\beta}$ in \cite{BD2009}. 
It is immediate that Dirac measures do not belong to ${\cal M}_\beta^+$. 
However absolutely continuous measures with respect to the Lebesgue measure, with density $\varphi\in L^2(\real^d)$ with compact support, do belong to ${\cal M}_\beta^+$ and will play an important role in the small-balls scaling. 
In that case, we shall abusively write $\mu\in L_c^2(\real^d)$ (here, the index $c$ stand for compact support). 
Recall the following properties on ${\cal M}_\beta^+$ from Propositions 2.2 and 2.3 from~\cite{BD2009}: 
\begin{prop}
\label{prop:Mcal}
\begin{enumerate}[(i)]
\item The set ${\cal M}_\beta^+$ is an affine subspace and, for all $\mu\in {\cal M}_\beta^+$,
\begin{equation*}
\label{maj:chap3}
 \int_{\real^d\times\real_+} \mu\big(B(x,r)\big)^2 r^{-\beta-1}\ dxdr <+\infty.
\end{equation*} 
\item If $d<\beta<2d$, then $L_c^2(\real^d)\subset {\cal M}_\beta^+$ and for all $\mu\in L_c^2(\real^d)$:
\begin{equation*}
\label{prop:L1L2}
\int_{\real^d}\mu\big(B(x,r)\big)^2\ dx\leq C_\mu\big( r^{d}\wedge r^{2d}\big).
\end{equation*}
\end{enumerate}
\end{prop}
Moreover, ${\cal M}_{\beta}$ is closed under translations, rotations and dilatations and is included in the subspace of diffuse measures, see Proposition~2.3 and Proposition~2.4 in \cite{BD2009} for details. 
See also \cite[Sec. 2.2]{KLNS2007} for a sufficient condition to belong to ${\cal M}_{\beta}$ in terms of the Riesz energy of a measure.

%%%%%%%%%%%%%%%%%%%%%%%%%%%%%%%%%%%%

\subsection*{Poissonian asymptotics}
Since our strategy consists in identifying in our functional Poissonian terms to which well known asymptotics are applied, recall these Poissonian asymptotics, from \cite{KLNS2007} (but with our current notations, see also \cite{BEK2010, BD2009}). 
\begin{theo}[Poissonian asymptotics, \cite{BEK2010}, \cite{BD2009} or \cite{KLNS2007}]
\label{theo:Poisson}
Let $\Phi$ be a marked Ppp in \eqref{eq:procM} and \eqref{eq:procMM} with compensator $K(0)dxF(dr)$ with $F$ having density $f$ satisfying  \eqref{eq:tails} for $d<\beta<2d$.

\begin{enumerate}[(i)]

\item{Large-balls scaling}: Assume $\lambda(\rho)\rho^\beta\to +\infty$. 
Then, for $n(\rho)=\big(\lambda(\rho)\rho^\beta\big)^{1/2}$, $\widetilde M_\rho(\cdot)/n(\rho)$ converges in the fdd sense on ${\cal M}_\beta^+$ to $W$  where
$$
W(\mu)=\int_{\real^d\times\real_+} \mu\big(B(x,r)\big)\ M_2(dx,dr)
$$
and $M_2$ is a centered Gaussian random measure with control measure $K(0)C_\beta r^{-\beta-1}\ dxdr$.

\item{Intermediate scaling}: Assume $\lambda(\rho)\rho^\beta\to a^{d-\beta}\in]0,+\infty[$. 
Then, for $n(\rho)=1$, $\widetilde M_\rho(\cdot)/n(\rho)$ converges in the fdd sense on ${\cal M}_\beta^+$ to $\widetilde P\circ D_a$ where  
$$
\widetilde P(\mu)=\int_{\real^d\times\real_+} \mu\big(B(x,r)\big)\ \widetilde \Pi(dx,dr)
$$
with $\widetilde \Pi$ a (compensated) Ppp with compensator measure $K(0)C_\beta r^{-\beta-1}\ dxdr$ and $D_a$ is the dilatation defined by $(D_a\mu)(B)=\mu(a^{-1}B)$.

\item{Small-balls scaling}: Assume $\lambda(\rho)\rho^\beta\to 0$. 
Then, for $n(\rho)=\big(\lambda(\rho)\rho^\beta\big)^{1/\gamma}$ with  $\gamma=\beta/d\in]1,2[$, $\widetilde M_\rho(\cdot)/n(\rho)$ converges in the fdd sense in $L^1(\real^d)\cap L^2(\real^d)$ to $Z$ where
$$
Z(\mu)=\int_{\real^d} \phi(x)\ M_\gamma(dx) \quad \mbox{ for } \quad\mu(dx)=\phi(x)dx \mbox{ with } \phi\in L^1(\real^d)\cap L^2(\real^d),
$$
is a stable integral with respect to a $\gamma$-stable random measure $\sigma_\gamma dx$ for 
$$
\sigma_\gamma=\frac{K(0) C_\beta v_d^\gamma}d \int_{0}^{+\infty} \frac{1-\cos(r)}{r^{1+\gamma}}dr
$$
and with unit skewness.
\end{enumerate}
\end{theo}
Here, and in the sequel, we follow the notations of \cite{ST1994} for stable random variables and integrals.

%%%%%%%%%%%%%%%%%%%%%%%%%%%%%%%%%%%%

\subsection{Large-balls scaling}
\label{sec:large_scaling}
In this section, we first investigate the behavior of $\widetilde{M}_\rho^R(\mu)$ in \eqref{eq:procMR} under the large-balls scaling $\lambda(\rho)\rho^\beta\to +\infty$ when $\rho \to 0$. 
As explained previously, the superposition due to $\lambda\to+\infty$ acts firstly producing a Gaussian field $W^R$ with a CLT type argument. 
Next, let $R\to +\infty$ to obtain the asymptotic behavior of $\widetilde{M}_\rho(\mu)$ according to \eqref{eq:intervertion}. 
The field obtained is given by a Gaussian integral similar to that of Theorem~\ref{theo:Poisson} (see also Theorem~2 (i) in \cite{KLNS2007}).

\begin{theo}[Large-balls scaling asymptotics]
\label{theo:largeBalls}
Assume \eqref{eq:tails} and the kernel $K_\rho$ satisfies  \eqref{eq:lambdarho}, \eqref{eq:controlintKrho} and Hypothesis~\ref{hyp:det1} for its associated operator ${\bf K}_\rho$ in \eqref{eq:KK}.
Suppose $\lambda(\rho)\rho^\beta\to +\infty$ when $\rho\to 0$, then the field $n(\rho)^{-1}\widetilde{M}_\rho(\cdot)$ converges in finite-dimensional distributions sense to $W(\cdot)$ in the space ${\cal M}_\beta^+$ where .
$$
W(\mu)=\int_{\real^d\times\real_+} \mu\big(B(x,r)\big)\ M_2(dx,dr),
$$
with a centered Gaussian random measure $M_2$ with control measure $K(0)C_\beta r^{-\beta-1}\ dxdr$.
\end{theo}

\medskip\noindent
Following our strategy, we start with the asymptotics of $\widetilde{M}_\rho^R(\mu)$:
\begin{prop}
\label{prop:cvmrlb}
Suppose $\lambda(\rho)\rho^\beta\to +\infty$ when $\rho \to 0$ and let $n(\rho)=\big(\lambda(\rho)\rho^\beta\big)^{1/2}$. 
Then, for all fixed $R>0$ and for all $\mu \in {\cal M}_\beta^+$, $n(\rho)^{-1}\widetilde{M}_\rho^R(\mu)$ converges in distribution when $\rho \to 0$ to
$$
W^R(\mu)=\int_{\real^d\times\real_+} \mu\big(B(x,r)\big)\ind_{\{r\leq R\}}\ M_2(dx,dr),
$$
uniformly in $R$, where $M_2$ is the same centered Gaussian random measure as in Theorem~\ref{theo:largeBalls}. 
\end{prop}
\begin{Proof}
The convergence in distribution of $\widetilde{M}_\rho^R(\mu)$ for $\mu\in{\cal M}_\beta^+$ is shown by the convergence of its Laplace transform: for $\theta\in\real$
\begin{equation}
\label{eq:Lapflulb}
\ee\Big[\exp\big(-\theta n(\rho)^{-1}\widetilde{M}_\rho^R(\mu)\big)\Big]
=\exp\big(\theta\ee[n(\rho)^{-1}M_\rho^R(\mu)]\big)\
\ee\Big[\exp\big(-\theta n(\rho)^{-1}M_\rho^R(\mu)\big)\Big].
\end{equation} 
Since $M_\rho^R$ given in \eqref{eq:procMR} is a determinantal integral with a compactly supported (say in $\Lambda_\mu^R$) integrand 
$g_\mu^R(x,r):=\mu\big(B(x,r)\big)\ind_{\{r\leq R\}}$, and the kernel $K$ satisfying Hypothesis~\ref{hyp:det1}, its Laplace transform is given by Theorem~\ref{theo:Laplace_Dpp}:
\begin{eqnarray} 
\nonumber
\ee\Big[\exp\big(-\theta n(\rho)^{-1}M_\rho^R(\mu)\big)\Big]
&=&\Det\Big( I-\widehat{K}_\rho\big[1-e^{-\theta n(\rho)^{-1} g_\mu^R}\big]\Big)\\
\label{eq:LT1lb}
&=&\exp\bigg(-\sum_{n\geq 1} \frac{1}{n}\Tr\Big( \widehat{K}_\rho\left[1-e^{-\theta n(\rho)^{-1}g_\mu^R}\right]^n\Big)\bigg),
\end{eqnarray}
where $\widehat{K}_\rho\big[1-e^{-\theta n(\rho)^{-1}g_\mu^R}\big]$ is the bounded operator of $L^2(\real^d\times \real_+)$ given in \eqref{eq:Kcrochet}. 
We compute the first trace in the sum in \eqref{eq:LT1lb} with Proposition~\ref{prop:int-tr} applied with the Dpp $\Phi_\rho$  with kernel $\widehat{K}_\rho$ on $\real^d\times\real_+$ restricted on the compact $\Lambda_\mu^R$ and the function $1-e^{-\theta n(\rho)^{-1}g_\mu^R}$ (see Proposition~\ref{prop:int-tr}):
\begin{eqnarray*} 
\Tr\Big(\widehat{K}_\rho\Big[1-e^{-\theta n(\rho)^{-1}g_\mu^R}\Big]\Big)
&=&\ee\left[\int_{\real^d\times \real_+}(1-e^{-\theta n(\rho)^{-1}g_\mu^R})\ \Phi_\rho(dx,dr)\right]\\
&=&\int_{\real^d\times \real_+}\big(1-e^{-\theta n(\rho)^{-1}\mu(B(x,r))\ind_{\{r\leq R\}}}\big)
K_\rho(x,x)\frac{f(r/\rho)}{\rho}\ dxdr.
\end{eqnarray*}
With \eqref{eq:lambdarho}, this term for $n=1$ combines with the factor  $\exp\big(\theta\ee[n(\rho)^{-1}M_\rho^R(\mu)]\big)$ of \eqref{eq:Lapflulb} into 
$$
\exp\bigg(\int_{\real^d\times\real_+}\psi\big(\theta n(\rho)^{-1}g_\mu^R\big)\ \lambda(\rho) K(0)\frac{f(r/\rho)}{\rho}\ dxdr\bigg)
$$
with $\psi(u)=e^{-u}-1+u$. 
The Laplace transform of $\widetilde{M}_\rho^R(\mu)$ in \eqref{eq:Lapflulb} thus rewrites
\begin{eqnarray} 
\nonumber
\ee\Big[\exp\Big(-\theta n(\rho)^{-1}\widetilde{M}_\rho^R(\mu)\Big)\Big]
\nonumber
&=&\exp\bigg(\int_{\real^d\times\real_+}\psi\big(\theta n(\rho)^{-1}g_\mu^R\big)\ \lambda(\rho) K(0)\frac{f(r/\rho)}{\rho}\ dxdr\bigg)\\
\label{eq:LT2lb}
&&\hspace{.8cm} \times\exp\bigg(-\sum_{n\geq 2} \frac{1}{n}\Tr\Big( \widehat{K}_\rho\Big[1-e^{-\theta n(\rho)^{-1}g_\mu^R}\Big]^n\Big) \bigg).
\end{eqnarray}
First, we deal with the first exponential term in \eqref{eq:LT2lb}: the key point is that this is the Laplace transform of $n(\rho)^{-1}\widetilde{P}_\rho^R(\mu)$ with
\begin{equation}
\label{eq:Prflulb}
\widetilde{P}_\rho^R(\mu)=\int_{\real^d\times\real_+}\mu\big(B(x,r)\big)\ind_{\{r\leq R\}}\ \widetilde{\Pi}_\rho(dx,dr),
\end{equation}
where $\widetilde{\Pi}_\rho$ is a compensated Poisson random measure on $\real^d\times\real_+$ with intensity 
$$
\lambda(\rho) K(0)\frac{f(r/\rho)}{\rho}\ dxdr.
$$ 
From (i) in Theorem~\ref{theo:Poisson} (Theorem~2-(i) in \cite{KLNS2007}), \eqref{eq:Prflulb} converges in distribution  when $\rho \to 0$ to the Gaussian integral $W^R(\mu)$. 
We show now that this convergence is actually uniform in $R$, to that way, consider the difference of the $\log$-Laplace transform of $n(\rho)^{-1}\widetilde{P}_\rho^R(\mu)$ and of $W^R(\mu)$:
\begin{eqnarray}
\nonumber
\lefteqn{\Big| \log\Big(\ee\Big[ \exp\big(n(\rho)^{-1}\widetilde{P}_\rho^R(\mu)\Big]\Big)-
\log\Big(\ee\Big[ \exp\big(W^R(\mu)\big)\Big]\Big)\Big|
}
\\
\nonumber
&\leq&\bigg|\int_{\real^d\times\real_+} \psi\big(n(\rho)^{-1}\mu\big(B(x,r)\big)\ind_{\{r\leq R\}}\big)\lambda(\rho) K(0)\frac{f(r/\rho)}{\rho}-\frac{\mu\big(B(x,r)\big)^2}{2}\ind_{\{r\leq R\}} \frac{K(0)}{ r^{\beta +1}} \ dxdr\bigg| 
\\
\nonumber
&\leq&\int_{\real^d\times\real_+} \Big| \psi\big(n(\rho)^{-1}\mu\big(B(x,r)\big)\ind_{\{r\leq R\}}\big)\lambda(\rho) K(0)\frac{f(r/\rho)}{\rho}-\frac{\mu\big(B(x,r)\big)^2}{2}\ind_{\{r\leq R\}} \frac{K(0)}{ r^{\beta +1}} \Big|\ dxdr. 
\\\label{eq:bound3}
\end{eqnarray}
Since $\psi(u)\sim \frac{u^2}{2}$ when $u\to 0$ and since $n(\rho)=\big(\lambda(\rho)\rho^\beta\big)^{1/2} \to +\infty$ when $\rho \to 0$, using the tails behaviour~\eqref{eq:tails}, we have:
\begin{eqnarray*}
\psi\Big(n(\rho)^{-1}\mu\big(B(x,r)\big)\Big)\lambda(\rho) K(0)\frac{f(r/\rho)}{\rho}&\underset{\rho \to 0}{\sim} &
\frac{\mu\big(B(x,r)\big)^2}{2n(\rho)^2}\lambda(\rho)K(0)\frac{\rho^\beta}{ r^{\beta+1}} \\
&&=\frac{\mu\big(B(x,r)\big)^2}{2} K(0)\frac{1}{ r^{\beta +1}},
\end{eqnarray*}
proving that the integrand in \eqref{eq:bound3} converges to $0$. 
Moreover, using \eqref{eq:tails} and $\psi(x)\leq x^2/2$, for all $r$ and for all $\rho>0$, we have:
\begin{eqnarray*}
\lefteqn{
\bigg|\psi\Big(n(\rho)^{-1}\mu\big(B(x,r)\big)\Big)\lambda(\rho) K(0)\frac{f(r/\rho)}{\rho}-\frac{\mu\big(B(x,r)\big)^2}{2} \frac{K(0)}{ r^{\beta +1}}\bigg|
}\\
&\leq& \frac{\mu\big(B(x,r)\big)^2}{2n(\rho)^2}\lambda(\rho) K(0)\frac{f(r/\rho)}{\rho}+ \frac{\mu\big(B(x,r)\big)^2}{2} \frac{K(0)}{ r^{\beta +1}} \notag \\
&\leq& K(0)(C_0+1)\frac{\mu\big(B(x,r)\big)^2}{2 r^{\beta+1}}
\end{eqnarray*}
which is integrable over $\real^d\times\real_+$ according to Proposition \ref{prop:Mcal}. 
Then, the dominated convergence theorem ensures that \eqref{eq:bound3} converges to $0$ when $\rho \to 0$. 
Moreover, since it does not depend on $R$, the convergence of $n(\rho)^{-1}\widetilde{P}_\rho^R(\mu)$ to $W^R(\mu)$ is uniform in $R$.

\bigskip\noindent
Next, we deal with the other second exponential terms in \eqref{eq:LT2lb} and show that they converge to $1$ proving that for all $n\geq 2$,
$$
\lim_{\rho \to 0} \Tr\Big(\widehat{K}_\rho\big[1-e^{-\theta n(\rho)^{-1}g_\mu^R}\big]^n\Big)=0.
$$
More precisely, the convergence of \eqref{eq:LT2lb} to $1$ will derive from the following lemmas. 
Recall $g_\mu^R(x,r)=\mu\big(B(x,r)\big)\ind_{\{r\leq R\}}$ and $\mu\in{\cal M}_\beta^+$; 
in particular $g_\mu^R$ is bounded with compact support. 
Since $K_\rho$ satisfies Hypothesis~\ref{hyp:det1},
Proposition~\ref{prop:hyp} first ensures  $\widehat{K}_\rho$ satisfies also Hypothesis~\ref{hyp:det1} and Proposition~\ref{prop:HS} next ensures that $\widehat{K}_\rho\big[1-e^{-\theta n(\rho)^{-1}g_\mu^R}\big]$ is the kernel of an Hilbert-Schmidt operator in \eqref{eq:KK}. 
\begin{lemme}
\label{lemme:trace2n}
For all $n\geq 2$, we have 
\begin{equation*}
\label{eq:Tr2n}
\Tr\Big( \widehat{K}_\rho\big[1-e^{-\theta n(\rho)^{-1}g_\mu^R}\big]^n\Big)
\leq \Tr\Big(\widehat{K}_\rho\big[1-e^{-\theta n(\rho)^{-1}g_\mu^R}\big]^2\Big)^{n/2}.
\end{equation*}
\end{lemme}
\begin{lemme}
\label{lemme:trace2}
Assume Conditions \eqref{eq:cfvol} and \eqref{eq:controlintKrho}, and consider $\mu \in{\cal M}_\beta^+$.
Then there is $\rho^*>0$ and a constant $C_K\in ]0,+\infty[$ such that for all $\rho\in ]0,\rho^*[$, uniformly in $R$,  
$$
\Tr\Big( \widehat{K}_\rho\big[1-e^{-\theta n(\rho)^{-1}g_\mu^R}\big]^2\Big)\leq
 C_K C_\mu C_f\theta^2\frac{\lambda(\rho)\rho^q}{n(\rho)^2}  
$$
with $C_f=\big(\int_0^{+\infty}  r^{q/2} f(r) dr \big)^2$.
\end{lemme}
As a consequence of both Lemma~\ref{lemme:trace2n}, we have
\begin{eqnarray}
\nonumber
\bigg| -\sum_{n\geq 2} \frac{1}{n}\Tr\Big(\widehat{K}_\rho\big[1-e^{-\theta n(\rho)^{-1}g_\mu^R}\big]^n\Big)\bigg| 
\nonumber
&\leq& \sum_{n\geq 1} \frac 1n\bigg(\sqrt{\Tr\Big( \widehat{K}_\rho\big[1-e^{-\theta n(\rho)^{-1}g_\mu^R}\big]^2\Big)}\bigg)^n \\
\label{eq:bound_log2}
&=&-\ln\bigg(1-\sqrt{\Tr\Big( \widehat{K}_\rho\big[1-e^{-\theta n(\rho)^{-1}g_\mu^R}\big]^2\Big)}\bigg).
\end{eqnarray}
Next, since \eqref{eq:cfvol} holds true under \eqref{eq:tails}, Lemma~\ref{lemme:trace2} applies and the bound \eqref{eq:bound_log2} goes to $0$ when $\rho \to 0$ since $\lambda(\rho)\rho^q/n(\rho)^2=\rho^{q-\beta}$ with $q>\beta$. 
As a consequence, 
$$
\lim_{\rho\to 0} \exp\bigg(-\sum_{n\geq 2} \frac 1n \Tr\Big( \widehat{K}_\rho\big[1-e^{-\theta n(\rho)^{-1}g_\mu^R}\big]^n\Big) \bigg)=1,
$$
and the limit in \eqref{eq:LT2lb} writes
\begin{equation*}
\label{eq:cvrholb}
\lim_{\rho \to 0} 
\ee\Big[\exp\big(-\theta n(\rho)^{-1}\widetilde{M}_\rho^R(\mu)\big)\Big]
=\ee\Big[\exp\big(-\theta W^R(\mu)\big)\Big],
\end{equation*}
achieving the proof of Proposition~\ref{prop:cvmrlb}. 
\CQFD
\end{Proof}

\medskip\noindent
It remains to prove Lemma~\ref{lemme:trace2n} and Lemma~\ref{lemme:trace2}.

\medskip\noindent
\begin{Proof} (Lemma~\ref{lemme:trace2n})
Recall that for a Hilbert-Schmidt operator $T$ with operator norm  $\|T\|$ and  Hilbert-Schmidt norm $\|T\|_2$, we have $\|T\| \leq \|T\|_2$ (see for instance Lemma~2.1 in \cite{ShiT2003} or \cite{RS1972} for details). 
Then, we have
\begin{eqnarray*}
\Tr\Big(\widehat{K}_\rho\big[1-e^{-\theta n(\rho)^{-1}g_\mu^R}\big]^n\Big) &\leq& \Big\|\widehat{K}_\rho\big[1-e^{-\theta n(\rho)^{-1}g_\mu^R}\big]\Big\|^{n-2}  \Tr\Big(\widehat{K}_\rho\big[1-e^{-\theta n(\rho)^{-1}g_\mu^R}\big]^2\Big)
\\
&\leq& \Big\|\widehat{K}_\rho\big[1-e^{-\theta n(\rho)^{-1}g_\mu^R}\big] \Big\|_2^{n-2}  \Tr\Big(\widehat{K}_\rho\big[1-e^{-\theta n(\rho)^{-1}g_\mu^R}\big]^2\Big).
\end{eqnarray*} 
Moreover we have:
\begin{eqnarray*}
\Big\|\widehat{K}_\rho\big[1-e^{-\theta n(\rho)^{-1}g_\mu^R}\big] \Big\|_2^2
&=&\int_{(\real^d\times\real_+)^{2}} \Big|\widehat{K}_\rho\big[1-e^{-\theta n(\rho)^{-1}g_\mu^R}\big] \big((x,r),(y,s)\big) \Big|^2 \ dxdydrds \\
&=&\int_{(\real^d\times\real_+)^{2}} \big(1-e^{-\theta n(\rho)^{-1}g_\mu^R(x,r)}\big) \widehat{K}_\rho\big((x,r),(y,s)\big)^2  \\
&&\hspace{4cm} \times\big(1-e^{-\theta n(\rho)^{-1}g_\mu^R(y,s)}\big)\ dxdydrds\\
&=& \Tr\Big(\widehat{K}_\rho\big[1-e^{-\theta n(\rho)^{-1}g_\mu^R}\big]^2\Big).
\end{eqnarray*}
and thus, we obtain, for every $n \geq 2$:
$$
\Tr\Big( \widehat{K}_\rho\big[1-e^{-\theta n(\rho)^{-1}g_\mu^R}\big]^n\Big)
\leq \Tr\Big(\widehat{K}_\rho\big[1-e^{-\theta n(\rho)^{-1}g_\mu^R}\big]^2\Big)^{n/2}.
$$
\CQFD
\end{Proof}

\medskip\noindent
\begin{Proof}(Lemma~\ref{lemme:trace2})
The operator $\widehat{K}_\rho\big[1-e^{-\theta n(\rho)^{-1}g_\mu^R}\big]^2$ is an integral operator with kernel
\begin{eqnarray*}
\lefteqn{
\widehat{K}_\rho\big[1-e^{-\theta n(\rho)^{-1} g_\mu^R}\big]^2 \big((x,r),(y,s)\big)}\\
&=&
\sqrt{1-e^{-\theta n(\rho)^{-1}g_\mu^R(x,r)}}\
\sqrt{\frac{f(r/\rho)}{\rho}}
\left(\int_{\real^d\times\real_+}\hskip -20pt \big(1-e^{-\theta n(\rho)^{-1}g_\mu^R(z,t)}\big)\frac{f(t/\rho)}{\rho} K_\rho(x,z) K_\rho(z,y)\ dzdt\right)
\\
&& \hspace{7cm}\times\sqrt{1-e^{-\theta n(\rho)^{-1}g_\mu^R(y,s)}}\ \sqrt{\frac{f(s/\rho)}{\rho}}.
\end{eqnarray*}
Its trace is thus given by:
\begin{eqnarray}
\nonumber
\Tr\Big( \widehat{K}_\rho\big[1-e^{-\theta n(\rho)^{-1}g_\mu^R}\big]^2\Big)
&=&\int_{\real^d\times\real_+} \widehat{K}_\rho\big[1-e^{-\theta n(\rho)^{-1}g_\mu^R}\big]^2\big((x,r),(x,r)\big)\ dxdr 
\\
\nonumber
&=&\int_{(\real^d\times\real_+)^2} \hskip -10pt \big(1-e^{-\theta n(\rho)^{-1}g_\mu^R(x,r)}\big) \big(1-e^{-\theta n(\rho)^{-1}g_\mu^R(z,t)}\big) \\
\label{eq:tr2lb}
&&\hspace{2cm} \times\frac{f(r/\rho)}{\rho}\frac{f(t/\rho)}{\rho} | K_\rho(x,z)|^2 \ dxdzdrdt.
\end{eqnarray} 
Since $\mu$ has a compact support, the function $g_\mu^R$ has also a compact support and 
$g_\mu^R(x,r)=0$ for, say, $\|r\|\geq M$. 
Thus the integrand in \eqref{eq:tr2lb} is a positive function with compact support (for $\theta$ or $\rho$ small enough). 
Dealing first with the integral over $\real^d \times \real^d$, 
since $1-e^{-\theta n(\rho)^{-1}g_\mu^R(x,r)}\leq \theta n(\rho)^{-1}\mu\big(B(x,r)\big)\ind_{\{r\leq R\}}$, we have 
\begin{eqnarray}
\nonumber 
\lefteqn{\int_{\real^d \times \real^d} \big(1-e^{-\theta n(\rho)^{-1}g_\mu^R(x,r)}\big) \big(1-e^{-\theta n(\rho)^{-1}g_\mu^R(z,t)}\big) |K_\rho(x,z)|^2\ dxdz
} \\
\nonumber
&=&\int_{B(0,M)\times B(0,M)} \big(1-e^{-\theta n(\rho)^{-1}g_\mu^R(x,r)}\big) \big(1-e^{-\theta n(\rho)^{-1}g_\mu^R(z,t)}\big) |K_\rho(x,z)|^2\ dxdz\\
\nonumber
&\leq & 
\int_{B(0,M)\times B(0,M)} \Big(\frac{\theta}{ n(\rho)}\Big)^2 \mu\big(B(x,r)\big)\mu\big(B(z,t)\big)\ind_{\{r\leq R\}}\ind_{\{t\leq R\}} K_\rho(x-z)^2\ dxdz\\
\nonumber
&\leq& \frac{\theta^2}{n(\rho)^2}\ind_{\{r\leq R\}}\ind_{\{t\leq R\}} \Big(\int_{B(0,M)\times B(0,M)}  \mu\big(B(x,r)\big)^2 K_\rho(x-z)^2\ dxdz\Big)^{1/2}\\
&&
\label{eq:CauSchlb1}
\hspace{4.5cm}\times \Big(\int_{B(0,M)\times B(0,M)}  \mu\big(B(z,t)\big)^2K_\rho(x-z)^2\ dxdz\Big)^{1/2}
\end{eqnarray}
using the Cauchy-Schwarz inequality. 
But, with the Fubini theorem, we have
\begin{eqnarray}
\nonumber
\int_{\real^d\times\real^d}  \mu\big(B(x,r)\big)^2 K_\rho(x-z)^2\ dxdz
&\leq&\int_{\real^d}  \mu\big(B(x,r)\big)^2  \Big(\int_{\real^d}K_\rho(x-z)^2\ dz\Big)\ dx\\
\label{eq:tech_CS11}
&\leq& C_K\lambda(\rho) C_\mu(r^p\wedge r^q)
\end{eqnarray}
since $\mu\in {\cal M}_\beta^+$ and 
using condition \eqref{eq:controlintKrho}.
Plugging into \eqref{eq:CauSchlb1}, \eqref{eq:tech_CS11} and a similar bound for the second integral in \eqref{eq:CauSchlb1}, we have 
\begin{eqnarray*}
\lefteqn{\int_{\real^d \times \real^d} \big(1-e^{-\theta n(\rho)^{-1}g_\mu^R(x,r)}\big) \big(1-e^{-\theta n(\rho)^{-1}g_\mu^R(z,t)}\big) |K_\rho(x,z)|^2\ dxdz
} \\
\nonumber
&\leq&C_K\theta^2\frac{\lambda(\rho)}{n(\rho)^2}\ind_{\{r\leq R\}}\ind_{\{t\leq R\}}C_\mu r^{q/2} t^{q/2}.
\end{eqnarray*}
As a consequence, the bound \eqref{eq:tr2lb} continues as follows
\begin{eqnarray*}
\nonumber
\lefteqn{\Tr\Big( \widehat{K}_\rho\big[1-e^{-\theta n(\rho)^{-1}g_\mu^R}\big]^2\Big)}
\\
&\leq &C_K\theta^2\frac{\lambda(\rho)}{n(\rho)^2}\int_{(\real^+)^2}\ind_{\{r\leq R\}}\ind_{\{t\leq R\}}C_\mu r^{q/2} t^{q/2} \frac{f(r/\rho)}{\rho}\frac{f(t/\rho)}{\rho} \ dr dt\\
&=&C_KC_\mu\theta^2\frac{\lambda(\rho)}{n(\rho)^2} \Big(\int_0^R r^{q/2}\frac{f(r/\rho)}{\rho}\ dr \Big)^2
=C_KC_\mu\theta^2\frac{\lambda(\rho)\rho^q}{n(\rho)^2} \Big(\int_0^{R/\rho}  r^{q/2} f(r) dr \Big)^2.
\end{eqnarray*}
But since $f$ is integrable and $q\leq 2d$ (Definition~\ref{def:Mbeta}) the finite volume condition \eqref{eq:cfvol} entails
$$
\int_0^{R/\rho} r^{q/2} f(r)\ dr\leq C_f:=\int_0^{+\infty} r^{q/2} f(r)\ dr <+\infty.
$$
\CQFD
\end{Proof}

\bigskip\noindent
We continue following the strategy exposed page~\pageref{sec:strategy}. 
Since the convergence in $\rho$ in Proposition~\ref{prop:cvmrlb} is uniform in $R$, the interchange \eqref{eq:intervertion} applies and we obtain:
\begin{equation*}
\label{intervertlb}
\lim_{\rho \to 0} {\cal L}\big(n(\rho)^{-1}\widetilde{M}_\rho(\mu)\big)
=\lim_{\rho \to 0}\lim_{R\to +\infty}  {\cal L}\big(n(\rho)^{-1}\widetilde{M}_\rho^R(\mu)\big)
=\lim_{R\to +\infty} {\cal L}\big(W^R(\mu)\big).
\end{equation*}
It remains now to identify $\lim_{R\to +\infty} W^R(\mu)$, this is done in the following proposition:
\begin{prop}
\label{prop:cvWR}
For all $\mu \in {\cal M}_\beta^+$, $W^R(\mu)$ converges in probability when $R \to +\infty$ to
$$
W(\mu)=\int_{\real^d\times\real_+} \mu\big(B(x,r)\big)\ M_2(dx,dr)
$$
where $M_2$ is the same centered Gaussian random measure as in Theorem~\ref{theo:largeBalls}. 
\end{prop}
\begin{Proof}
Since $W^R(\mu)$ and $W(\mu)$ are both integral with respect to the same Gaussian measure $M_2$, we have:
$$
W(\mu)-W^R(\mu)=\int_{\real^d\times\real_+} \mu\big(B(x,r)\big) \ind_{\{ r> R\}}\ M_2(dx,dr)
$$
whose $\log$-Laplace transform is
\begin{equation}
\label{eq:LT2011W}
\log\Big(\ee\Big[\exp\big(W(\mu)-W^R(\mu)\big)\Big]\Big)
=\frac{1}{2}\int_{\real^d\times\real_+} \mu\big(B(x,r)\big)^2\ \ind_{\{r> R\}} K(0)r^{-\beta-1}\ dxdr.
\end{equation}
The integrand in \eqref{eq:LT2011W} converges to $0$ when $R\to +\infty$ and is bounded by 
$$
\mu\big(B(x,r)\big)^2 K(0) r^{-\beta-1}
$$ 
which, thanks to Proposition~\ref{prop:Mcal}, is integrable for $\mu \in {\cal M}_\beta^+$. 
The dominated convergence theorem thus ensures that \eqref{eq:LT2011W} converges to $0$, i.e. $W(\mu)-W^R(\mu)\Rightarrow 0$ and $W^R(\mu)\stackrel{\pp}{\longrightarrow} W(\mu)$, $R\to+\infty$, which is Proposition~\ref{prop:cvWR}.
\CQFD
\end{Proof}

\bigskip\noindent
So far, all the intermediate results are obtained to prove Theorem~\ref{theo:largeBalls}:

\medskip\noindent
\begin{Proof}[Th. \ref{theo:largeBalls}]
The one-dimensional convergence is obtained by the combination of \eqref{eq:intervertion} with Proposition~\ref{prop:cvmR}, Proposition~\ref{prop:cvmrlb} and Proposition~\ref{prop:cvWR}. 
Now, remark that the fields $\widetilde{M}_\rho$ and $W$ are both linear on ${\cal M}_\beta^+$. 
Thus, using the Cram\'er-Wold device and the linear structure of $\mathcal{M}_\beta$, we have immediately the convergence of the finite-dimensional distributions from the one-dimensional convergence.
\CQFD
\end{Proof}

%%%%%%%%%%%%%%%%%%%%%%%%%%%%%%%%%%%%%%%%%%%%%%%%%%%%

\subsection{Intermediate scaling}
\label{sec:intermediate_scaling}

This section investigates the asymptotic behavior of $\widetilde M_\rho$  in \eqref{eq:procMM} under the intermediate scaling, when $\lim_{\rho\to+\infty}\lambda(\rho) \rho^\beta= a\in]0,+\infty[$. 
In this section, set $n(\rho)=1$.

\begin{theo}[Intermediate scaling asymptotics]
\label{theo:intermediateBalls}
Assume \eqref{eq:tails} and the kernel $K_\rho$ satisfies  \eqref{eq:lambdarho}, \eqref{eq:controlintKrho} and Hypothesis~\ref{hyp:det1} for its associated operator ${\bf K}_\rho$ in \eqref{eq:KK}.
Suppose $\lambda(\rho)\rho^\beta\to a^{d-\beta} \in ]0, +\infty[$ when $\rho\to 0$, then $\widetilde{M}_\rho(\cdot)$ converges in the finite-dimensional distributions sense to $\widetilde{P}\circ D_a(\cdot)$ in the space ${\cal M}_\beta^+$, where 
$$
\widetilde{P}(\mu)=\int_{\real^d\times\real_+}\mu\big(B(x,r)\big)\ \widetilde{\Pi}(dx,dr),
$$
with $\widetilde{\Pi}$ a compensated Poisson random measure on $\real^d\times\real_+$ with intensity measure $K(0)C_\beta r^{-\beta-1}dxdr$ and $D_a$ standing for the dilatation $(D_a\mu)(B)=\mu(a^{-1}B)$. 
\end{theo}

\medskip\noindent
Following the same strategy as previously (see page \pageref{sec:strategy}), first investigate the asymptotic behavior of $\widetilde{M}_\rho^R(\mu)$ in \eqref{eq:procMR} when $\rho\to 0$ and next let $R\to +\infty$ in the obtained limit. 
Roughly speaking, as in the Poissonian case (see (ii) in Theorem~\ref{theo:Poisson}, or Theorem~2-(ii) in \cite{KLNS2007}), the limit corresponds to take the limit in the intensity of the underlying random measure.
The result states as follows
\begin{prop}
\label{prop:cvmri}
Suppose $\lambda(\rho)\rho^\beta \to a \in (0,+\infty)$ when $\rho \to 0$. 
Then, for all $\mu \in {\cal M}_\beta^+$ and $R>0$, $\widetilde{M}_\rho^R(\mu)$ converges in distribution to 
$$
\big(\widetilde{P}^R\circ D_a\big)(\mu)=\int_{\real^d\times\real_+}\mu\big(B(x,r)\big)\ind_{\{r\leq R\}}\ \widetilde{\Pi}(dx,dr),
$$
where $\widetilde{\Pi}$ is the same compensated Poisson random measure as in Theorem~\ref{theo:intermediateBalls}.
\end{prop}
\begin{Proof}
The proof follows the same scheme as for Proposition \ref{prop:cvmrlb}.
Recall that in this context, $n(\rho)=1$ is set. 
The Laplace transform of $\widetilde{M}_\rho^R(\mu)$ is given by \eqref{eq:LT2lb}, i.e. 
\begin{eqnarray}
\nonumber
\ee\Big[\exp\big(-\theta\widetilde{M}_\rho^R(\mu)\big)\Big]
&=&\exp\bigg(\int_{\real^d\times\real_+}\psi\big(\theta \mu\big(B(x,r)\big)\ind_{\{r\leq R\}}\big)K_\rho(x,x)\frac{f(r/\rho)}{\rho} dxdr\bigg)
\\
\label{eq:eq:LT2i}
&&\hspace{3cm} \times\exp\bigg(-\sum_{n\geq 2} \frac{1}{n}\Tr\big( \widehat{K}_\rho\big[1-e^{-\theta g_\mu^R}\big]^n\big)\bigg).\end{eqnarray}
The first exponential in \eqref{eq:eq:LT2i} is the Laplace transform of 
$$
\widetilde{P}_\rho^R(\mu)=\int_{\real^d\times\real_+}\mu\big(B(x,r)\big)\ind_{\{r\leq R\}}\ \widetilde{\Pi}_\rho(dx,dr),
$$
where $\widetilde{\Pi}_\rho$ is a compensated Poisson random measure on $\real^d\times\real_+$ with intensity measure $\lambda(\rho) K(0)\frac{f(r/\rho)}{\rho}dxdr$. 
From (ii) in Theorem~\ref{theo:Poisson} (see also Theorem~2-(i) in \cite{KLNS2007}), under Condition \eqref{eq:tails}, when $\lim_{\rho\to 0}\lambda(\rho)\rho^\beta = a^{d-\beta} \in ]0,+\infty[$, this process converges to
$$
\big(\widetilde{P}^R\circ D_a\big)(\mu)=\int_{\real^d\times\real_+}\mu\big(B(x,r)\big)\ind_{\{r\leq R\}}\ \widetilde{\Pi}(dx,dr),
$$
where $\widetilde{\Pi}$ is a compensated Poisson random measure on $\real^d\times\real_+$ with intensity measure $K(0)r^{-\beta-1}dxdr$.
In particular, we have : 
$$
\lim_{\rho\to 0}\exp\bigg(\int_{\real^d\times\real_+}\psi\Big(\theta \mu\big(B(x,r)\big)\ind_{\{r\leq R\}}\Big)K_\rho(x,x)\frac{f(r/\rho)}{\rho}\ dxdr\bigg)
=\ee\Big[\exp\big(-\theta\big(\widetilde{P}^R\circ D_a\big)(\mu)\big)\Big].
$$
The proof is completed by showing that the second exponential term in \eqref{eq:eq:LT2i} converges to $1$. 
Proceeding as in the proof of Proposition~\ref{prop:cvmrlb}, with $n(\rho)=1$, Lemma~\ref{lemme:trace2} entails
$$
\Tr\Big(\widehat{K}_\rho\big[1-e^{-\theta g_\mu^R}\big]^2\Big)
\leq C_K C_\mu C_f \theta^2 \lambda(\rho)\rho^q 
$$ 
which goes to $0$ since $\lim_{\rho\to 0}\lambda(\rho)\rho^q=0$ for $q>\beta$. 
As a consequence
$$
\lim_{\rho\to 0} \Tr\Big(\widehat{K}_\rho\big[1-e^{-\theta g_\mu^R}\big]^2\Big)=0.
$$ 
Then, with Lemma~\ref{lemme:trace2n}, we still have for every $n \geq 2$
$$
\Tr\Big(\widehat{K}_\rho\big[1-e^{-\theta g_\mu^R}\big]^n\Big)
\leq \Tr\Big( \widehat{K}_\rho\big[1-e^{-\theta g_\mu^R}\big]^2\Big)^{n/2},
$$
and the second exponential term in \eqref{eq:eq:LT2i} converges to $1$, as in the proof of Proposition~\ref{prop:cvmrlb}, page~\pageref{eq:cvrholb}, this concludes the proof of Proposition~\ref{prop:cvmri}.
\CQFD
\end{Proof}

\bigskip\noindent
Combining Proposition \ref{prop:cvmri} with the interchange \eqref{eq:intervertion}, we have:
$$
\lim_{\rho \to 0} {\cal L}\big(\widetilde{M}_\rho(\mu)\big)
=\lim_{R\to +\infty} \lim_{\rho \to 0} {\cal L}\big(\widetilde{M}_\rho^R(\mu)\big)
=\lim_{R\to +\infty}{\cal L}\big(\widetilde{P}^R(\mu)\big).
$$
It remains now to identify $\lim_{R\to +\infty}\widetilde{P}^R(\mu)$, this is done in the following proposition:
\begin{prop}
\label{prop:cvpRi}
For all $\mu\in {\cal M}_\beta^+$, $\widetilde{P}^R(\mu)$ converges in $L^1$ when $R\to +\infty$ to 
$$
\widetilde{P}(\mu)=\int_{\real^d\times\real_+}\mu\big(B(x,r)\big)\ \widetilde{\Pi}(dx,dr),
$$
where $\widetilde{\Pi}$ is the same compensated Poisson random measure as in Theorem~\ref{theo:intermediateBalls}.
\end{prop}
\begin{Proof}
Since $\widetilde{P}^R(\mu)$ and $\widetilde{P}(\mu)$ are Poissonian integral with respect to the same measure $\widetilde{\Pi}$, we have:
$$
\Big| \widetilde{P}^R(\mu)-\widetilde{P}(\mu)\Big|
=\bigg|\int_{\real^d\times\real_+} \mu\big(B(x,r)\big)\ind_{\{r> R\}} \ \widetilde{\Pi}(dx,dr)\bigg| 
$$
and
\begin{eqnarray*}
\ee\Big[\Big| \widetilde{P}^R(\mu)-\widetilde{P}(\mu)\Big|\Big]
&\leq& 2 \int_{\real^d\times\real_+} \mu\big(B(x,r)\big)\ind_{\{r> R\}} K(0)r^{-\beta-1}dxdr
\\
&\leq& 2v_d\mu(\real^d)K(0) \int_R^{+\infty} r^{d-\beta-1}dr
\\
&=&\frac{2v_d\mu(\real^d)K(0)}{(\beta-d)R^{\beta-d}}\longrightarrow 0, \: R\to +\infty. 
\end{eqnarray*}
\CQFD
\end{Proof}

\bigskip\noindent
So far, all the intermediate results are obtained to proove Theorem~\ref{theo:intermediateBalls}: 

\medskip\noindent
\begin{Proof}[Th.~\ref{theo:intermediateBalls}]
The one-dimensional convergence is obtained by the combination of \eqref{eq:intervertion} with Proposition \ref{prop:cvmR}, Proposition \ref{prop:cvmri} and Proposition \ref{prop:cvpRi}.
Since the fields $\widetilde{M}_\rho$ and $\widetilde{P}$ are both linear on ${\cal M}_\beta^+$, 
 using the Cram\'er-Wold device and the linear structure of $\mathcal{M}_\beta$, we have immediately the convergence of the finite-dimensional distributions from the one-dimensional convergence.
\CQFD
\end{Proof}

%%%%%%%%%%%%%%%%%%%%%%%%%%%%%%%%%%%%%%%%%%%%%%%%%%%%%%%%%%
 
\subsection{Small-balls scaling}
\label{sec:small_scaling}

This section investigates the asymptotics of $\widetilde{M}_\rho^R(\mu)$ under the small-balls scaling, i.e. when $\lim_{\rho\to 0} \lambda(\rho)\rho^\beta=0$. 
We deal first with the limit in $\rho$ of the truncated field $\widetilde{M}_\rho^R(\mu)$. 
In this case, the obtained limit does not depend on $R$, roughly speaking this is due to the fast decreasing of the rescaled radii $\rho r$ since $\rho\to 0$ very fast in this regime.  
The limiting field thus obtained is a stable integral similar to the one obtained for the Poissonian model in (iii) of Theorem~\ref{theo:Poisson} (cf. also Theorem~2-(iii) in \cite{KLNS2007} and cf. \cite{ST1994} for notations on stable integrals). 
In this case, the limit is driven by small balls and this requires to consider smooth configuration $\mu(dx)=\varphi(x)dx$. Roughly speaking, if the configuration $\mu$ were, for instance, atomic, there will be a possibility for the small balls driving the asymptotics to not charge $\mu$ and $M(\mu)$ would vanish. 

\begin{theo}
\label{theo:smallBalls}
Assume \eqref{eq:tails} and the kernel $K_\rho$ satisfies  \eqref{eq:lambdarho}, \eqref{eq:controlintKrho} and Hypothesis~\ref{hyp:det1} for its associated operator ${\bf K}_\rho$ in \eqref{eq:KK}.
Suppose $\lambda(\rho)\rho^\beta \to 0$ when $\rho \to 0$ and set $n(\rho)=(\lambda(\rho)\rho^\beta)^{1/\gamma}$ with $\gamma=\beta/d$. 
Then, the field $n(\rho)^{-1}\widetilde{M}_\rho(\cdot)$ converges in the finite-dimensional distributions sense when $\rho\to 0$ to $Z(\cdot)$ in $L_c^2(\real^d)$ where 
$$
Z(\mu)=\int_{\real^d}\varphi(x) M_\gamma(dx), \quad \mbox{ for } 
\mu(dx)=\varphi(x)dx
$$
with $M_\gamma$ a $\gamma$-stable measure with control measure $\sigma_\gamma dx$ where
\begin{equation*}
\label{eq:sigma}
\sigma_\gamma=\frac{K(0)C_\beta v_d^{\gamma}}{d}\int_0^{+\infty}\frac{1-\cos(r)}{r^{1+\gamma}}\ dr
\end{equation*}
and constant unit skewness.
\end{theo}

\medskip\noindent
First, we have:
\begin{prop}
\label{prop:cvmrsb}
Suppose $\lambda(\rho)\rho^\beta\to 0$ when $\rho \to 0$ and set $n(\rho)=\big(\lambda(\rho)\rho^\beta\big)^{1/\gamma}$ for $\gamma=\beta/d\in]1,2[$. 
Then, for all $R>0$ and for all $\mu \in L^1(\real^d)\cap L^2(\real^d)$, writing $\mu(dx)=\varphi(x)dx$, $n(\rho)^{-1}\widetilde{M}_\rho^R(\mu)$ converges in the finite-dimensional distributions sense when $\rho\to 0$ to
$$
Z(\mu)=\int_{\real^d}\varphi(x) M_\gamma(dx),
$$
where $M_\gamma$ is the same $\gamma$-stable measure as in Theorem~\ref{theo:smallBalls}.
\end{prop}
\begin{Proof}
Recall the Laplace transform of $\widetilde{M}_\rho^R(\mu)$ is given in \eqref{eq:LT2lb}:
\begin{eqnarray*}
\ee\Big[\exp\Big(-\theta n(\rho)^{-1}\widetilde{M}_\rho^R(\mu)\Big)\Big]
\nonumber
&=&\exp\bigg(\int_{\real^d\times\real_+}\psi\big(\theta n(\rho)^{-1}g_\mu^R\big)\ K_\rho(x,x)\frac{f(r/\rho)}{\rho}\ dxdr\bigg)\\
&&\hspace{.8cm} \times\exp\bigg(-\sum_{n\geq 2} \frac{1}{n}\Tr\Big( \widehat{K}_\rho\Big[1-e^{-\theta n(\rho)^{-1}g_\mu^R}\Big]^n\Big) \bigg).
\end{eqnarray*}
The first exponential term is still the Laplace transform of
$n(\rho)^{-1}\widetilde{P}_\rho(\mu)$ where $\widetilde{P}_\rho(\mu)$ is the compensated Poissonian integral \eqref{eq:Prflulb}. 
With the change of variable $r=n(\rho)^{1/d}s$, this $\log$-Laplace transform becomes:
\begin{equation}
\label{eq:LT2sb}
\int_{\real^d\times\real_+}\psi\Big(\theta n(\rho)^{-1}\mu\Big(B\big(x,n(\rho)^{1/d}s\big)\Big)\ind_{\{s<n(\rho)^{-1/d}R\}}\Big)\lambda(\rho)K(0)n(\rho)^{1/d}\frac{f\big(sn(\rho)^{1/d}/\rho\big)}{\rho}\ dxds. 
\end{equation}
For $\mu(dx)=\varphi(x)dx$ with $\varphi \in L_c^2(\real^d)$, then the following Lemma from \cite{KLNS2007} entails
$$
\lim_{\rho\to 0} \theta n(\rho)^{-1}\mu\Big(B\big(x,n(\rho)^{1/d}s\big)\Big)\ind_{\{ s<n(\rho)^{-1/d}R\}} =\theta\varphi(x)v_d s^d,
$$
$dx$-almost everywhere and 
$$
x\mapsto \sup_{r>0} \bigg(\frac{\mu\big(B(x,r)\big)}{v_d r^d}\bigg) \in L^2(\real^d).
$$

\begin{lemme}[Lemma 4 in \cite{KLNS2007}]
Let $C$ be a bounded Borelian set in $\real^d$ with $Leb(C)=1$. 
\begin{enumerate}[(i)]
\item If $\varphi\in L^1$, then 
$\lim_{v\to 0} v^{-1}\int_{x+v^{1/d} C}\varphi(y)\ dy=\varphi(x)$ for $dx$-almost all $x$.
\item If $\varphi\in L^1$, then 
$\varphi_*(x):=\sup_{v>0} v^{-1} \int_{x+v^{1/d}C}|\varphi(y)|\ dy<+\infty$ for $dx$-almost all $x$
\item Moreover if $\varphi\in L^p$ for some $p>1$ then $\varphi_*\in L^p$.
\end{enumerate}
\end{lemme}
Then, using the very argument of the proof of Theorem 2 in \cite{KLNS2007} (see also the proof of Theorem~2.16 in \cite{BD2009})
\begin{eqnarray}
\nonumber
\int_{\real^d\times\real_+}\psi\Big(\theta n(\rho)^{-1}\mu\big(B(x,n(\rho)^{1/d}r)\big)\ind_{\{r<n(\rho)^{-1/d}R\}}\Big )\lambda(\rho)K(0)n(\rho)^{1/d}\frac{f(r n(\rho)^{1/d}/\rho)}{\rho}\ dxdr
&&\\
\label{eqLT2sb}
\sim_{\rho \to 0}\lambda(\rho)K(0) \int_{\real^d\times\real_+}\psi\big(\theta\varphi(x)v_d r^d\big) n(\rho)^{1/d}\frac{f(r n(\rho)^{1/d}/\rho)}{\rho}\ dxdr.&&
\end{eqnarray}
Using now the proof of Theorem 2 in \cite{KLNS2007} under the small-ball scaling, the right-hand side in \eqref{eq:LT2sb} converges to the Laplace transform of $Z(\mu)$. 
This implies that the random variable $n(\rho)^{-1}\widetilde{P}_\rho(\mu)$ converges in distribution to $Z(\mu)$.

\medskip\noindent
The proof is completed by showing that the second exponential term in \eqref{eq:LT2sb} converges to $1$. 
Using the same conclusion as in the proof of Proposition~\ref{prop:cvmrlb} page~\pageref{eq:cvrholb} with Lemma~\ref{lemme:trace2n}, it is enough to show that for this regime we still have
$$
\lim_{\rho\to 0}\Tr\Big(\widehat{K}_\rho\big[1-e^{-\theta n(\rho)^{-1}g_\mu^R}\big]^2\Big)=0.
$$ 
Since we consider $\mu \in L_c^2(\real^d)$, we have also $\mu \in L^1(\real^d)$ and Proposition \ref{prop:Mcal}-(ii) ensures that we can take here $q=2d$ and then Lemma~\ref{lemme:trace2} writes with $n(\rho)=(\lambda(\rho)\rho^\beta)^{1/\gamma}$:
\begin{equation*}
\label{eq:majtr2sb}
\Tr\Big(\widehat{K}_\rho\big[1-e^{-\theta n(\rho)^{-1}g_\mu^R}\big]^2\Big)
\leq C_K C_\mu C_f\theta^2 \frac{\lambda(\rho)\rho^{2d}}{ n(\rho)^2}
=C_K C_\mu C_f\theta^2\lambda(\rho)^{(\beta-2d)/\beta}
\end{equation*}
which goes to $0$ when $\rho \to 0$ since $\beta<2d$. 
\CQFD
\end{Proof}

\bigskip\noindent
So far all the intermediate results are obtained to finish the proof of Theorem~\ref{theo:smallBalls} as for Theorem~\ref{theo:largeBalls} and Theorem~\ref{theo:intermediateBalls}.

%%%%%%%%%%%%%%%%%%%%%%%%%%%%%%%%%%%%%%%%%%%%%%%%%%%%%%%%%%%%%%%%%

\section{Examples}
\label{sec:Exemples}
In this section, we provide two examples of Dpps satisfying our hypotheses and illustrating our results of Section~\ref{sec:asymptotics}. 

%%%%%%%%%%%%%%%%%%%%%%%%%%%

\subsection{Ginibre process}
\label{sec:Ginibre}

The Ginibre point process $\phi$ is a Dpp with kernel
\begin{equation*}
\label{eq:ginker}
K^G(x,y)=\exp \Big(-\frac 12\|x-y\|^2\Big), \quad x,y\in \real^d,
\end{equation*}
with respect to the Lebesgue measure. 
Such processes have been used recently to model wireless networks of communication, see \cite{DZH}. 
For our macroscopic analysis, we introduce the following scaled version of the Ginibre point process. 
Let $\phi_\rho$ be a scaled Ginibre point process with kernel:
\begin{equation}
\label{eq:ginkerscal}
K_{\rho}^{G}(x,y)=\lambda (\rho) \exp \Big(-\frac{\lambda(\rho)^\delta}{2}\|x-y\|^2\Big), \quad x,y\in \real^d,
\end{equation}
with respect to the Lebesgue measure, where $\lambda :\real_+ \rightarrow \real_+$ is a decreasing function with $\lim_{\rho\to 0}\lambda(\rho)=+\infty$ and with $ \delta \geq 2/d$ , so that \eqref{eq:controlintKrho} is satisfied. 
Then, Theorem~\ref{theo:largeBalls}, Theorem~\ref{theo:intermediateBalls} and Theorem~\ref{theo:smallBalls}  apply for the Ginibre Kernel \eqref{eq:ginkerscal} with $K(0)=1$ therein. 

\begin{Rem}
\label{rem:alphaGinibre}
{\rm 
The results apply also for the thinned and re-scaled Ginibre point process $\phi_\alpha$ (or $\alpha$-Ginibre point process) with kernel:
\begin{equation*}
\label{eq:ginkerthin}
K^{G,\alpha}(x,y)=\exp \Big(-\frac{\|x-y\|^2}{2\alpha}\Big),
\end{equation*}
where $0<\alpha\leq 1$. 
Such a process is obtained by retaining independently each point of the  Ginibre point process with probability $\alpha$ and then applying a scaling to conserve the density (mean number of points by volume unit) of the initial Ginibre point process.
This so-called $\alpha$-Ginibre point process bridges smoothly between the Ginibre point process ($\alpha=1$) and the Poisson point process ($\alpha\to 0$).  
For the scaled version, replace \eqref{eq:ginkerscal} by
\begin{equation*}
\label{eq:ginkerscalthin}
K_{\rho}^{G,\alpha}(x,y)=\lambda (\rho) \exp \Big(-\frac{\lambda(\rho)^\delta}{2\alpha}\|x-y\|^2\Big).
\end{equation*}
} 
\end{Rem}

%%%%%%%%%%%%%%%%%%%%%%%%%%%%%%%%%%%%%%%%%%%%%%%%%%

\subsection{Bessel process}
\label{sec:Bessel}

The Bessel process is a Dpp $\phi$ with kernel 
\begin{equation}
\label{eq:sincker}
K^B(x,y)=\frac{\sqrt{\Gamma (d/2+1)}}{\pi ^{d/4}}\frac{J_{d/2}\big(2\sqrt{\pi}\Gamma (d/2+1)^{1/d}\|x-y\|\big)}{\|x-y\|^{d/2}}, \quad x,y\in \real^d,
\end{equation}
with respect to the Lebesgue measure, where $J_{d/2}$ stands for the Bessel function of the first kind. 
For instance, for $d=1$ we have
$$
K^B(x,y)=\frac{\sin \big(\pi \|x -y\|\big)}{\pi \|x-y\|}.
$$
For our macroscopic analysis, we introduce the following scaled version of the Bessel point process. 
Let $\phi_\rho$ be a scaled Bessel point process with kernel:
\begin{equation}
\label{eq:sinckerscal}
K_{\rho}^B(x,y)=\frac{\sqrt{\lambda(\rho)\Gamma (d/2+1)}}{\pi ^{d/4}}\frac{J_{d/2}\big(2\sqrt{\pi}\Gamma (d/2+1)^{1/d}\lambda(\rho)^{1/d}\|x-y\|\big)}{\|x-y\|^{d/2}}
\end{equation}
with respect to the Lebesgue measure, where $\lambda :\real_+ \rightarrow \real_+$ is a decreasing function with $\lim_{\rho\to 0}\lambda(\rho)=+\infty$. 
Using the following asymptotics of the Bessel functions of the first kind (see \cite{Arfken_Weber}):
\begin{eqnarray*}
J_\alpha (r) &\underset{r\to 0}{\sim}&\frac{1}{\Gamma (\alpha +1)}\left(\frac{r}{2}\right)^\alpha,\\
J_\alpha (r) &\underset{r\to +\infty}{\sim}&\sqrt{\frac{2}{\pi r}}\cos \Big(r-\frac{\alpha \pi}{2}-\frac{\pi}{4} \Big),
\end{eqnarray*}
condition \eqref{eq:controlintKrho} is satisfied so that our main results Theorem~\ref{theo:largeBalls}, Theorem~\ref{theo:intermediateBalls} and Theorem~\ref{theo:smallBalls} apply (with $K(0)=1$).

%%%%%%%%%%%%%%%%%%%%%%%%%%%%%%%%%%%%%%%%%%%%%%%%%%%%%%%%%%%%%%%%%%%

\section{Comments}
\label{sec:comments}

\subsection{Zoom-in asymptotics}

For the Poisson random balls model, the study of the microscopic fluctuations in \cite{BE2006} obtained by zooming-in instead of zooming-out, leads to very similar results to those obtained in the macroscopic behavior in \cite{KLNS2007} under the large-ball scaling and the intermediate scaling. 
This similarity is the origin of the unified approach for both types of fluctuations in \cite{BEK2010}, used also in the weighted model in \cite{BD2009}. 
In the microscopic point of view, this is the behavior of small balls which matters and this is encapsulated in \cite{BEK2010} in the following condition on small radii 
$$
f(r)\sim_{r\to 0} \frac{1}{r^{\beta+1}}.
$$
In this case, $f$ cannot be a probability density nor be integrable. 
Consequently, we can not study a determinantal random balls model under a zoom-in procedure. 
Indeed, even if we were to consider a marked Dpp on $\real^d\times\real_+$ with kernel:
\begin{equation}
\label{eq:cfqueuezi}
\widehat{K}\big((x,r),(y,s)\big)=\sqrt{f(r)}K(x,y)\sqrt{f(s)}
\end{equation}
with a determinantal kernel $K$ on $\real^d$ and $f$ a function on $\real_+$ satisfying condition \eqref{eq:tails}, this Dpp would have no chance to satisfy Hypothesis~\ref{hyp:det1} when $f$ is not integrable.

%%%%%%%%%%%%%%%%%%%%%%%%%%%%%%%%%%%%%%%%%%%%%%%%%%

\subsection{$\alpha$-determinantal and $\alpha$-permanental processes}
\label{sec:alphaDet}

The Dpps actually belong to a larger class of point processes, the so-called $\alpha$-determinantal/per\-manental processes. 
When $\alpha>0$, such processes exhibit attraction between its particles, and when $\alpha<0$, they  exhibit repulsiveness. 
When $\alpha=-1$, the (usual) Dpp are recovered while the case $\alpha=1$ corresponds to permanental processes. 
The definition $\alpha$-determinantal/permanental processes follows the same definition as Def.~\ref{def:Dpp} but with the determinant replaced by a $\alpha$-determinant. 
Recall that for a matrix $A=(a_{i,j})_{1\leq i,j\leq n}$ and $\alpha\in\real$, its $\alpha$-determinant is defined by 
\begin{equation}
\label{eq:det_alpha}
{\det}_\alpha A
=\sum_{\sigma\in{\mathfrak S}_n}
\alpha^{n-\nu(\sigma)} \prod_{i=1}^n a_{i,\sigma(i)}
\end{equation}
where ${\mathfrak S}_n$ is the symmetric group of permutation of $\{1,\dots, n\}$ and $\nu(\sigma)$  is the number of cycles in $\sigma\in{\mathfrak S}_n$. 
When $\alpha=-1$ (resp. $\alpha=1$), \eqref{eq:det_alpha} defines the (standard) determinant (resp. permanent) of $A$ : ${\det}_{-1} A=\det A$, ${\det}_{1} A=\mbox{perm }A$. 

\medskip\noindent
The following result from \cite{ShiT2003} extends Theorem~\ref{theo:Laplace_Dpp} and proves the existence of such processes for some $\alpha$'s and gives their Laplace transform: 
\begin{theo}[Th.~1.2 in \cite{ShiT2003}]
\label{theo:Laplace_Det} 
Let $E$ be a Polish space equipped with a diffuse Radon measure $\lambda$ and $K$ be a bounded symmetric integral operator on $L^2(E,\lambda)$ satisfying Hypothesis~\ref{hyp:det1}. 
Then for $\alpha\in \{2/m : m\in \nit\}\cup \{-1/m : m\in\nit\}$, there exists a unique point process $\phi$ such that 
\begin{equation}
\label{eq:Laplace_adet}
\ee\Big[\exp\Big(-\int f(x)\ \phi(dx)\Big)\Big]
=\Det\Big(I+\alpha K\big[1-e^{-f}\big]\Big)
\end{equation}
for each compactly supported measurable $f:E\to\real_+$ where $K[1-e^{-f}]$ still stands for the kernel \eqref{eq:Kcrochet}. 
Moreover, $\phi$ is a simple point process whose correlation functions are given by 
$$
\rho_{n,\alpha, K}(x_1, x_2,\dots, x_n)
={\det}_\alpha\Big(\big(K(x_i,x_j)\big)_{1\leq i,j\leq n}\Big).
$$
\end{theo}
Like for \eqref{eq:Fredholm1}, for a trace-class operator $T$ with $\|\alpha T\|<1$, the Fredholm determinant of $I-\alpha T$ expands in terms of $\alpha$-determinant
$$
\Det\big(I-\alpha T\big)^{-1/\alpha}
=\sum_{n=0}^{+\infty} \int_{E^n}{\det}_\alpha\big( (T(x_i,x_j))_{1\leq i,j\leq n}\big)\ \lambda^{\otimes n}(dx_1, \dots, dx_n). 
$$

\medskip\noindent
Using the expansion \eqref{eq:Fredholm_det} of the Fredholm determinant of the Laplace transform \eqref{eq:Laplace_adet}, our arguments can be carried out similarly for $\alpha$-determinantal/permanental processes. 
Indeed, since $|\alpha|\leq 1$, the terms for $n\geq 2$ can be similarly bounded and are still asymptotically negligible while the term $n=1$ is obviously the same Poissonian term. 
As a consequence, Theorems~\ref{theo:largeBalls}, \ref{theo:intermediateBalls}, \ref{theo:smallBalls} have natural generalization to $\alpha$-determinantal/permanental processes.

%%%%%%%%%%%%%%%%%%%%%%%%%%%%%%%%%%%%%%%%%%%%%%%%%%%%%%%%%%%%%%%%%%%%
%
\subsection{Non-stationary determinantal random ball model}

With slight modifications, our main results remain true for non-stationary determinantal random ball models.
Consider a determinantal process $\phi$ with kernel $K(x,y)$ still satisfying Hypothesis~\ref{hyp:det1} but also 
\begin{eqnarray}
\label{eq:Kxx}
x \longmapsto K(x,x) \in L^{\infty}(\real^d)
\end{eqnarray}
The zoom-out procedure consists now in introducing the scaled version $\phi_\rho$ of $\phi$, with kernel $K_\rho$ with respect to the Lebesgue measure satisfying
$$
K_\rho(x,x)\underset{\rho \rightarrow 0}{\sim}\lambda (\rho)K(x,x)
$$
with $\lim_{\rho\to 0} \lambda (\rho)= +\infty$. 
We also replace \eqref{eq:lambdarho} and \eqref{eq:controlintKrho} 
by 
\begin{eqnarray}
\label{eq:lambdarho2}
\sup_{x\in\real^d} K_\rho(x,x)\leq \lambda (\rho)\sup_{x\in\real^d} K(x,x)\\
\label{eq:controlintKrho2} 
\sup_{x\in\real^d}\int_{\real^d}\big|K_\rho(x,y) \big|^2 \ dy
\underset{\rho \rightarrow 0}{=}\mathcal{O}\big(\lambda(\rho)\big).
\end{eqnarray}
In this non-stationary context, Theorem~\ref{theo:largeBalls}, \ref{theo:intermediateBalls}, \ref{theo:smallBalls} have the following counterparts: 

\begin{theo} 

Assume \eqref{eq:tails} and $\phi_\rho$ is a Dpp with kernel satisfying \eqref{eq:Kxx}, \eqref{eq:lambdarho2}, \eqref{eq:controlintKrho2} and Hypothesis~\ref{hyp:det1} for its associated operator ${\bf K}_\rho$ in \eqref{eq:KK}.

\begin{enumerate}[(i)]

\item{Large-balls scaling}: Assume $\lambda(\rho)\rho^\beta\to +\infty$. 
Then, for $n(\rho)=\big(\lambda(\rho)\rho^\beta\big)^{1/2}$, $\widetilde M_\rho(\cdot)/n(\rho)$ converges in the fdd sense on ${\cal M}_\beta^+$ to $W$  where
$$
W(\mu)=\int_{\real^d\times\real_+} \mu\big(B(x,r)\big)\ M_2(dx,dr)
$$
and $M_2$ is a centered Gaussian random measure with control measure $K(x,x)C_\beta r^{-\beta-1}\ dxdr$.

\item{Intermediate scaling}: Assume $\lambda(\rho)\rho^\beta\to a^{d-\beta}\in]0,+\infty[$. 
Then, for $n(\rho)=1$, $\widetilde M_\rho(\cdot)/n(\rho)$ converges in the fdd sense on ${\cal M}_\beta^+$ to $\widetilde P\circ D_a$ where  
$$
\widetilde P(\mu)=\int_{\real^d\times\real_+} \mu\big(B(x,r)\big)\ \widetilde \Pi(dx,dr)
$$
with $\widetilde \Pi$ a (compensated) Ppp with compensator measure $K(x,x)C_\beta r^{-\beta-1}\ dxdr$ and $D_a$ is the dilatation defined by $(D_a\mu)(B)=\mu(a^{-1}B)$.

\item{Small-balls scaling}: Suppose $\lambda(\rho)\rho^\beta \to 0$ when $\rho \to 0$ and set $n(\rho)=(\lambda(\rho)\rho^\beta)^{1/\gamma}$ with $\gamma=\beta/d$. 
Then, the field $n(\rho)^{-1}\widetilde{M}_\rho(\cdot)$ converges in the finite-dimensional distributions sense when $\rho\to 0$ to $Z(\cdot)$ in $L_c^2(\real^d)$ where 
$$
Z(\mu)=\int_{\real^d}\varphi(x) M_\gamma(dx), \quad \mbox{ for } 
\mu(dx)=\varphi(x)dx
$$
with $M_\gamma$ a $\gamma$-stable measure with control measure $\sigma_\gamma K(x,x) dx$ where
\begin{equation*}
\label{eq:sigma}
\sigma_\gamma=\frac{C_\beta v_d^{\gamma}}{d}\int_0^{+\infty}\frac{1-\cos(r)}{r^{1+\gamma}}\ dr
\end{equation*}
and constant unit skewness.
\end{enumerate}
\end{theo}

\medskip\noindent
In this non-stationary case, the proof follows the same general strategy as in page~\pageref{sec:strategy} but with technical details requiring  \eqref{eq:Kxx}, \eqref{eq:lambdarho2}, \eqref{eq:controlintKrho2}. 
Roughly speaking, the limits are driven by the term $n=1$ in \eqref{eq:LT1lb} while the other terms ($n\geq 2$) are still negligible. 
Note that, in this non-stationary setting,  the Poissonian limits for $n=1$ come now from \cite{Gobard2014} (with $G=\delta_1$ therein) instead of \cite{KLNS2007}. 
Details are left to the interested readers.

%%%%%%%%%%%%%%%%%%%%%%%%%%%%%%%%%%%%%%%%%%%%%%%%%%%%%%%%%%%%%%%%%%%%

\appendix

\section{Appendix: (Marked) Determinantal Point Processes} 
\label{sec:mDpp}

In this section, we give a short presentation of Determinantal point processes (Dpps). 
For a general reference on point processes, refer to the two volumes book \cite{DVJ1} and for a specific reference on Dpps, refer to \cite{HKPV2009}. 
Dpps form a special class of point processes that exhibit repulsiveness between its points. 
Below, we consider a point process $\xi$ that is a (random) collection of locally finite points in, say, some Polish space $E$. 
In the sequel, to avoid any ambiguity, the points of the process are called particles.
In the following, simple point processes for which almost surely its points are all distinct are considered.
Considering a reference Borel measure $\mu$ on $E$, the distribution law of $\xi$ is, in general, characterized by its joint intensities.
\begin{df}
\label{def:intensity}
The joint intensities of a point process $\xi$ on a Polish space $E$ with respect to $\mu$ are functions (if any exists) $\rho_k : E \rightarrow [0,+\infty[$, $k\geq 1$, such that for any family of mutually disjoint Borelian subsets $D_1, \dots, D_k$ of $E$,
$$
\ee\bigg[\prod_{i=1}^k \xi(D_i) \bigg]
=\int_{\prod_{i=1}^k D_i}\rho_k(x_1, \dots,x_k)\ \mu(dx_1)\dots\mu(dx_k).
$$
\end{df}
Roughly speaking, $\rho_k(x_1, \dots,x_k)$ can be interpreted as the (infinitesimal) probability for $\xi$ to have particles in each $x_1, \dots, x_k$. 
For a Poisson point process (Ppp), the intensity are constant.
For a Dpp, the joint intensities are given by a certain determinant of a kernel $K$ characterizing the process, hence its name. 
\begin{df}
\label{def:Dpp}
A point process $\xi$ on $E$ is said to be a determinantal point process with kernel $K$ if it is simple and its joint intensities write for all $k\geq 1$ and all $ x_1, \dots, x_k \in E$:
\begin{equation*}
\label{eq:Dpp}
\rho_k(x_1, \dots, x_k)=\det\big(K(x_i,x_j)\big)_{1\leq i,j\leq k}:=\det[K](x_1, \dots, x_k).
\end{equation*}
\end{df}
See below in Theorem~\ref{theo:Laplace_Dpp} for condition ensuring the existence of such processes. 
Observe that the repulsiveness exhibited by a Dpp can be read on its joint intensity of second order. 
Indeed, if $K$ is continuous and  $x_1, x_2\in E$, the more they will be close to each other, the more the determinant of $\big(K(x_i,x_j) \big)_{1\leq i,j\leq 2}$ will be close to $0$. 
Thus, $\rho_2(x_1, x_2)\approx 0$ whenever $x_1\approx x_2$. 
This implies that, if there is a particle of the process in $x_1$, the probability that there is another particle in the close vicinity of $x_1$ is small. 
For a Ppp, the constant intensities show that the particles of a Dpp are independently drawn. 

\medskip\noindent
An important class of Dpp is the class of those whose kernel satisfies special properties (see Hypothesis~\ref{hyp:det2} below) that we recall from \cite{HKPV2009} and in our setting is encapsulated in Hypothesis~\ref{hyp:det1}. 
For that purpose, recall that, for a compact operator $T$ on a separable Hilbert space $H$ equipped with the scalar product $\langle\cdot,\cdot\rangle$, its trace is given by 
$$
\Tr(T)=\sum_{n=1}^{+\infty} \langle Te_n,e_n\rangle
$$
where $(e_n)_{n\geq 1}$  is (any) complete orthonormal (CONB) system of $H$. 
In particular, $T$ is said to be a trace-class operator if 
$$
\|T\|_1:=\Tr\big(|T|\big)<+\infty
$$
where $|T|=\sqrt{T^*T}$. 
The hypothesis on the kernel $K$ writes: 
\begin{hyp} 
\label{hyp:det2}
The Polish space $E$ is equipped with a Radon $\sigma$-finite measure $\lambda$. 
The map ${\bf K}$ is an operator from $L^2(E,\lambda)$ into $L^2(E,\lambda)$ satisfying the following conditions:
\begin{enumerate}[(i)]
\item \label{Hyp1}
${\bf K}$ is a bounded symmetric integral operator on $L^2(E,\lambda)$ with kernel $K$, i.e., for any $x\in E$ and any $f\in L^2(E,\lambda)$, 
$$
{\bf K}f(x)=\int_{E}K(x,y)f(y)\ \lambda(dy).
$$
\item \label{Hyp2}
 The spectrum of ${\bf K}$ is included in $[0,1[$. 

\item \label{Hyp3}
The map ${\bf K}$ is locally trace-class, i.e. for all compact $\Lambda \subset E$, the restriction ${\bf K}_\Lambda$ of ${\bf K}$ on $L^2(\Lambda,\lambda)$ is of trace-class. 
\end{enumerate}
\end{hyp}

\begin{Rem}
{\rm
If $K$ is the kernel of a map ${\bf K}$ satisfying Hypothesis~\ref{hyp:det2}, then $x\mapsto K(x,x)$ is nonnegative.
}
\end{Rem}
n the sequel, the limit in distribution of quantities \eqref{eq:Mrhomu1} is investigated by considering the Laplace transform of a Dpp. 
It is given in Theorem~\ref{theo:Laplace_Dpp} above from \cite{ShiT2003} and expressed in terms of Fredholm determinant. 
Recall that if $\|T\|_1<1$, the Fredholm determinant of $I+T$ is given by 
\begin{equation}
\label{eq:Fredholm_det}
\Det(I+T)=\exp\bigg(\sum_{n=1}^{+\infty} \frac{(-1)^{n-1}}n\Tr\big(T^n\big)\bigg).
\end{equation}
Moreover, for a trace-class operator $T$ with $\|T\|<1$, the Fredholm determinant of $I+T$ expands in terms of determinants as follows
\begin{equation}
\label{eq:Fredholm1}
\Det\big(I+T\big)
=\sum_{n=0}^{+\infty} \int_{E^n}{\det}\big( (T(x_i,x_j))_{1\leq i,j\leq n}\big)\ \lambda^{\otimes n}(dx_1, \dots, dx_n). 
\end{equation}

\begin{theo}[Th.~1.2 in \cite{ShiT2003}]
\label{theo:Laplace_Dpp}
Let $E$ be a Polish space equipped with a diffuse Radon measure $\lambda$ and $K$ be a bounded symmetric integral operator on $L^2(E,\lambda)$ satisfying Hypothesis~\ref{hyp:det2}. 
Then there exists a unique Dpp $\phi$ as in Definition \ref{def:Dpp} and its Laplace transform is given for each compactly supported measurable $f:E\to\real_+$ by 
\begin{equation}
\label{eq:Laplace_Dpp}
\ee\Big[\exp\Big(-\int f(x)\ \phi(dx)\Big)\Big]
=\Det\Big(I-K\big[1-e^{-f}\big]\Big)
\end{equation}
where $K[1-e^{-f}]$ stands for the kernel
\begin{equation}
\label{eq:Kcrochet}
K\big[1-e^{-f}\big](x,y)=\sqrt{1-\exp(-f(x))}\ K(x,y)\ \sqrt{1-\exp(-f(y))}.
\end{equation}
\end{theo}

\medskip\noindent
The following result is obtained by differentiation of the Laplace transform:
\begin{prop}
\label{prop:int-tr}
Let $\phi$ be a Dpp on a Polish space $E$ with kernel $K$ satisfying Hypothesis \ref{hyp:det2} with respect to a measure $\lambda$ on $E$. 
For any compact set $\Lambda$ of $E$ and any non-negative function $f$ defined on $E$, we have
$$
\ee\Big[\int_{\Lambda}f d\phi\Big]
=\int_\Lambda f(x) K(x,x)\ \lambda(dx)
=\Tr \big(K_\Lambda[f]\big).
$$
\end{prop}

\bigskip\noindent
In Section~\ref{sec:DetBall}, marked determinantal point processes are considered and, for that purpose some useful results on marked Dpps are given in the rest of this section.  
First, the following classical result on Ppps (see for instance Lemma 6.4.VI in \cite{DVJ1}) is easily extended:
{\it If $\xi=\{X_i\}_{i\geq 1}$ is a Ppp on a Polish space $E$ with intensity $\lambda\in\real_+$ and $(R_i)_{i\geq 1}$ is a family of {\it iid} random variables with distribution $F$ on a Polish space $E'$ (independent of $\xi$), then $\xi'=(X_i,R_i)_{i\geq 1}$ is a Ppp on $E\times E'$ with intensity $\lambda\otimes F$}.
In the determinantal case, we have:
\begin{prop}
\label{prop:MDPP}
Let $\phi=(X_i)_{i\geq 1}$ be a determinantal point process on a Polish space $E$ with kernel $K$, with respect to a Radon measure $\lambda$, and let $(R_i)_{i\geq 1}$ be a family of {\it iid} random variables on $\real_+$, independent of $(X_i)_{i\geq 1}$, with probability density $f$. 
Let $\Phi=\big\{(X_i,R_i)\big\}_{i\geq 1}$. 
Then, $\Phi$ is a determinantal point process on $E\times\real_+$ with kernel 
\begin{equation}
\label{eq:Kchap}
\widehat{K}\big((x,r),(y,s)\big)=\sqrt{f(r)}K(x,y)\sqrt{f(s)},
\end{equation}
with respect to the measure $\lambda(dx)dr$.
\end{prop}
The result still holds true for marks with values in a Polish space but in the sequel,only positive marks are used (i.e. $R_i\in\real_+$). 

\medskip\noindent
\begin{Proof}
To prove that $\Phi$ is a Dpp with kernel $\widehat{K}$, the joint intensities are shown to write
\begin{equation*}
\hat{\rho}_n\big((x_1,r_1),\dots,(x_n,r_n)\big)=\det\Big(\widehat{K}\big((x_i,r_i),(x_j,r_j)\big)_{1\leq i,j\leq n}\Big).
\end{equation*}
For all $n\geq 1$ and all set $A$, the symbol 
$\sum_{a_1,\dots,a_n \in A}^{\neq}$ 
will stand for the sum over all $n$-tuples $(a_1,\dots,a_n)\in A$ with pairwise distinct $a_i$ ($a_i\neq a_j$ for $i\neq j$ in $\{1,\dots,n\}$). 
Let $n\geq 1$ and $h$ a Borel function from $(E\times \real_+)^n$ to $\real_+$. 
We have :
\begin{eqnarray*}
\lefteqn{
\ee\bigg[ \sum_{(x_1,r_1),\dots,(x_n,r_n) \in \Phi}^{\neq}  h\big((x_1,r_1),\dots,(x_n,r_n)\big)\bigg]
}\\
&=&\ee\bigg[\sum_{x_1,\dots,x_n \in \phi}^{\neq} h\big((x_1,R_1),\dots,(x_n,R_n)\big)\bigg]\\
&=&\ee\left[ \ee\bigg[\sum_{x_1,\dots,x_n \in \phi}^{\neq} h\big((x_1,R_1),\dots,(x_{i_n},R_{i_n})\big)\Big\vert \phi \bigg] \right] \\
&=&\ee\left[ \int_{(\real_+)^n}\sum_{x_1,\dots,x_n \in \phi}^{\neq} h\big((x_1,r_1),\dots,(x_n,r_n)\big)\ \prod_{1\leq i\leq n} f(r_i)dr_i\right] \\
&=&\ee\left[ \sum_{x_1,\dots,x_n \in \phi}^{\neq} \int_{(\real_+)^n} h\big((x_1,r_1),\dots,(x_n,r_n)\big)\ \prod_{1\leq i\leq n} f(r_i)dr_i\right] \\
&=&\int_{(E\times \real_+)^n}  h\big((x_1,r_1),\dots,(x_n,r_n)\big)\prod_{1\leq i\leq n} f(r_i) \rho_n(x_1,\dots,x_n)\lambda(dx_1)dr_1\dots\lambda(dx_n) dr_n,
\end{eqnarray*}
where $\rho_n(x_1,\dots,x_n)=\Det[K](x_1,\dots,x_n)$ is the joint intensity of order $n$ of the Dpp $\phi$. 
Now, note that 
$$
\prod_{1\leq i\leq n} f(r_i)\ \Det[K](x_1,\dots,x_n)
=\Det\big[\widehat{K}\big]\big((x_1,r_1),\dots,(x_n,r_n)\big),
$$
where $\widehat{K}$ is given in \eqref{eq:Kchap}. 
Then
\begin{eqnarray*}
\lefteqn{
\ee\bigg[ \sum_{(x_1,r_1),\dots,(x_n,r_n) \in \Phi}^{\neq}  h\big((x_1,r_1),\dots,(x_n,r_n)\big)\bigg]
}\\
&=&
\int_{(E\times\real_+)^n} h\big((x_1,r_1),\dots,(x_n,r_n)\big) \ \Det\big[\widehat{K}\big]\big((x_1,r_1),\dots,(x_n,r_n)\big)\ \prod_{1\leq i\leq n} \lambda(dx_i)dr_i ,
\end{eqnarray*}
and, according to Definition~\ref{def:intensity} and Definition~\ref{def:Dpp}, $\Phi$ is a Dpp on $E\times\real_+$ with kernel $\widehat{K}$ with respect to the measure $\lambda(dx)dr$.
\CQFD
\end{Proof}

\medskip\noindent
Next, in the case where ${\bf K}$ satisfies Hypothesis~\ref{hyp:det2}, the operator $\widehat{{\bf K}}$ associated to $\widehat{K}$ defined in \eqref{eq:Kchap} above inherits these properties:
\begin{prop}
\label{prop:hyp}
Let ${\bf K}$ be an operator on $L^2(E,\lambda)$ satisfying Hypothesis~\ref{hyp:det2} and $\widehat{{\bf K}}$ be the integral operator with kernel \eqref{eq:Kchap} with probability density $f$.
Then, $\widehat{{\bf K}}$ satisfies Hypothesis~\ref{hyp:det2}.
\end{prop}
\begin{Proof}
We show that each point of Hypothesis~\ref{hyp:det2} is satisfied.
\begin{enumerate}[(i)]
\item $\widehat{{\bf K}}$ is obviously a symmetric integral operator and it is bounded since it is an Hilbert-Schmidt operator.
\item Let $\gamma\in [0,1[$ be in the spectrum of $\widehat{{\bf K}}$ and $g_\gamma$ an associated eigenfunction. 
Then,
\begin{eqnarray*}
\gamma g_\gamma(x,r)&=&\widehat{{\bf K}}g_\gamma(x,r)\\
&=&\int_{E\times \real_+} \sqrt{f(r)}K(x,y)\sqrt{f(s)}\ g_\gamma(y,s) \ \lambda(dy)ds \\
&=&\sqrt{f(r)}\int_E K(x,y) \int_{\real_+} \sqrt{f(s)} g_\gamma(y,s)\ ds \lambda(dy) \\
&=&\sqrt{f(r)} K\left(\int_{\real_+} \sqrt{f(s)}g_{\gamma}(\cdot,s)ds\right)(x).
\end{eqnarray*}
Thus, since $f$ is a probability density, 
\begin{eqnarray*}
\gamma \int_{\real_+} \sqrt{f(r)} g_\gamma(x,r)dr&=&\int_{\real_+} f(r) K\left(\int_{\real_+} \sqrt{f(s)} g_{\gamma}(\cdot,s)ds\right)(x)\ dr \\
&=&\int_{\real_+} f(r)dr\ K\left(\int_{\real_+} \sqrt{f(s)}g_{\gamma}(\cdot,s)ds\right)(x) \\
&=&K\left(\int_{\real_+} \sqrt{f(s)}g_{\gamma}(\cdot,s)ds\right)(x),
\end{eqnarray*}
proving that $\gamma$ is in the spectrum of $K$ (associated to the eigenfunction \\
 $x \mapsto \int_{\real_+}  \sqrt{f(r)} g_\gamma(x,r)dr$) and obviously $\gamma \in [0,1[$.
 
 \item First, let $\Lambda=\Lambda_E\times \Lambda_{\real_+}$ be a compact of $E\times \real_+$ and $\widehat{\bf K}_\Lambda$ be the restriction of $\widehat{\bf K}$ on $\Lambda$. 
In order to compute the trace of $\widehat{\bf K}_\Lambda$, consider a complete orthogonal basis (CONB in short) of $L^2\big(\Lambda, \lambda(dx)dr\big)$. 
Let $(e_n)_{n\geq 1}$, resp. $(b_n)_{n\geq 1}$, be a CONB of $L^2(\Lambda_E,\lambda)$, resp. of $L^2(\Lambda_{\real_+},dr)$. 
Then $(h_{n,k})_{n,k\geq 1}$,with $h_{n,k}(x,r)=e_n(x)b_k(r)$ is a CONB of $L^2\big(\Lambda, \lambda(dx)dr\big)$ (see \cite{RS1972}) and
$$
\Tr\big(\widehat{\bf K}_\Lambda\big)=\sum_{n,k\geq 1} \big\langle\widehat{\bf K}_\Lambda h_{n,k}, h_{n,k}\big\rangle_{L^2(\Lambda, \lambda(dx)dr)},
$$
 with for $n,k \geq 1$:
 \begin{eqnarray*}
 \lefteqn{\big\langle\widehat{\bf K}_\Lambda h_{n,k}, h_{n,k}\big\rangle_{L^2(\Lambda, \lambda(dx)dr)}
 }\\
&=&\int_{\Lambda^2}h_{n,k}(x,r) \widehat{\bf K}_\Lambda h_{n,k}(x,r) \ \lambda(dx)dr\\
&=&\int_{\Lambda^2}e_n(x)b_k(r)\sqrt{f(r)}K(x,y)\sqrt{f(s)}e_n(y)b_k(s)
\ \lambda(dy)ds \lambda(dx)dr\\
 &=&\left(\int_{\Lambda_{\real_+}}\sqrt{f(r)}b_k(r)dr \right)^2 \left(\int_{\Lambda_E^2}e_n(x)K(x,y)e_n(y)\ \lambda(dx)\lambda(dy)\right)\\
 &\leq& 
\langle \sqrt{f}, b_k\rangle_{L^2(\real_+)}^2 \big\langle Ke_n,e_n\big\rangle_{L^2(\Lambda_E)}
 \end{eqnarray*}
with the Fubini theorem. 
As a consequence, with the Bessel inequality, 
$\Tr\big(\widehat{\bf K}_\Lambda\big)\leq \|\sqrt{f}\|_{L^2(\real_+)}^2\Tr\big({\bf K}_{\Lambda_E}\big) <+\infty$,
 and $\widehat{\bf K}_\Lambda$ is locally trace-class. 
 Note that it is still true for subset $\Lambda$ of the form $\Lambda_E \times \real_+$.
\\
Next, for a general compact set $\Lambda$ of $E\times \real_+$, we have $\Lambda\subset \Lambda_E\times \Lambda_{\real_+}$ for compact sets $\Lambda_E$ of $E$ and $\Lambda_{\real_+}$ of $\real_+$.
Using the reunion $(m_n)_{n\geq 1}=(c_n)_{n\geq 1}\cup(d_n)_{n\geq 1}$ of orthonormal basis  $(c_n)_{n\geq1}$ of $L^2\big(\Lambda, \lambda(dx)dr\big)$ and $(d_n)_{n\geq 1}$ of $L^2\big(\Lambda_E\times \Lambda_{\real_+}\setminus\Lambda, \lambda(dx)dr\big)$, we have an orthnormal basis of  $L^2\big(\Lambda_E\times \Lambda_{\real_+},\lambda(dx)dr\big)$ and by the first part
\begin{eqnarray*}
\Tr\big(\widehat{\bf K}_{\Lambda_E \times \Lambda_{\real_+}}\big) &=&\sum_{n\geq 1} \big\langle\widehat{\bf K} m_n, m_n\big\rangle_{L^2(\Lambda_E \times \Lambda_{\real_+}, \lambda(dx)dr)}
\\
&=&
\sum_{n\geq 1} \big\langle\widehat{\bf K} c_n, c_n\big\rangle_{L^2(\Lambda, \lambda(dx)dr)}
+
\sum_{n\geq 1} \big\langle\widehat{\bf K} d_n, d_n\big\rangle_{L^2(\Lambda_E \times \Lambda_{\real_+}\setminus\Lambda, \lambda(dx)dr)}\\
&=&\Tr\big(\widehat{\bf K}_{\Lambda}\big)+
\Tr\big(\widehat{\bf K}_{\Lambda_E \times \Lambda_{\real_+}\setminus\Lambda}\big).
\end{eqnarray*}
Since all summands are positive we have $\Tr\big(\widehat{\bf K}_{\Lambda}\big)<+\infty$.
\end{enumerate}
\CQFD
\end{Proof}
\begin{Rem}
\label{rem:L1}
{\rm 
Straightforwardly, Proposition~\ref{prop:hyp} is still true for $f\in L^1(\real_+)$ but with condition \eqref{Hyp2} replaced by:
(ii') The spectrum of $\widehat{K}$ is included in $\big[0,\|f\|_1^{-1}\big[$. 
}
\end{Rem}

\begin{prop}
\label{prop:HS}
Let $K$ be a kernel satisfying Hypothesis~\ref{hyp:det2} and $g:E\to[0+\infty[$ be a bounded function with compact support. 
Then $K[g]$ given by
$$
K[g](x,y)=\sqrt{g(x)} K(x,y) \sqrt{g(y)}
$$
is the kernel of an Hilbert-Schmidt operator. 
\end{prop}
\begin{Proof}
The Hilbert-Schmidt property is shown by proving 
\begin{equation*}
\label{eq:Kg_HS}
\int_{E\times E} K[g](x,y)^2\ dxdy<+\infty.
\end{equation*}
Let $B$ be the compact support of $g$, using $\rho_2(x_1,x_2)=\det\big(K(x_i, x_j)_{1\leq i,j\leq 2}\big)\geq 0$, we have  
\begin{eqnarray*}
\int_{E\times E} K[g](x,y)^2\ \lambda(dx)\lambda(dy)
&=&\int_{E\times E} g(x) K(x,y)^2 g(y)\ \lambda(dx)\lambda(dy)\\
&\leq&\|g\|_\infty^2\int_{B\times B} K(x,x) K(y,y) \ \lambda(dx)\lambda(dy)\\
&=&\|g\|_\infty^2\Big(\int_{B} K(x,x) \ \lambda(dx)\Big)^2
\end{eqnarray*}
which is finite since $K$ is locally trace-class (Hypothesis~\ref{hyp:det2}).
\CQFD
\end{Proof}

%%%%%%%%%%%%%%%%%%%%%%%%%%%%%%%%%%%%%%%%%%%%%%%%%%%%%%%%%%%%%%%%%%%%

\section*{Acknowledgement}

The authors thank an anonymous referee for valuable comments on a preliminary draft of this manuscript.
The third author thanks also the \href{https://www.lebesgue.fr/}{Centre Henri Lebesgue} ($''$Investissements d'avenir$''$ program --- ANR-11-LABX-0020-01) for its financial support.

%%%%%%%%%%%%%%%%%%%%%%%%%%%%%%%%%%%%%%%%%%%%%%%%%%%%%%%%%%%%%%%%%%%%

%%%%%%%%%%%%%%%%%%%%%%%%%%%%%%%%%%%%%%%%%%%%%%%%%%%%%%%%%%%%%%%%%%%%%

\end{document}